%% file: _main.tex
\documentclass[reqno,letterpaper]{amsart}

\usepackage{thmtools}
\usepackage{thm-restate}
\usepackage{import, subfiles}
\usepackage{epsfig}
\usepackage{framed}
\usepackage{enumitem}

\usepackage{natbib}
\usepackage[linktoc=all,linkbordercolor={1 1 0}]{hyperref}
\usepackage{subcaption}
\usepackage{xspace}
\usepackage[normalem]{ulem}
\usepackage[usenames,dvipsnames]{color}
\usepackage{xr}
\usepackage{datetime}

\usepackage[capitalize]{cleveref}  

\crefname{lem}{Lemma}{Lemmas}
\crefname{prop}{Proposition}{Propositions}
\crefname{cor}{Corollary}{Corollaries}
\crefname{nthm}{Theorem}{Theorems}
\crefname{rmk}{Remark}{Remark}
\crefname{remark}{Remark}{Remark}
\crefname{assumption}{Assumption}{Assumptions}

\numberwithin{equation}{section}
\theoremstyle{plain}
\newtheorem{nthm}{Theorem}[section]
\theoremstyle{remark}

\newtheorem{assumption}{Assumption}
\makeatletter
\newcommand{\setassumptiontag}[1]{%
  \let\oldtheassumption\theassumption%
  \renewcommand{\theassumption}{#1}%
  \g@addto@macro\endassumption{%
    \addtocounter{assumption}{-1}%
    \global\let\theassumption\oldtheassumption}%
  }
\makeatother

\renewcommand*{\theassumption}{\arabic{section}.\Alph{assumption}}
\declaretheorem[style=remark,qed=$\bigtriangledown$,name=Remark,sibling=nthm]{remark}

\usepackage{misc}
\input{shortcuts2.tex}

\usepackage[usenames,dvipsnames]{color}
\usepackage[mathscr]{eucal} 
\usepackage{algorithmic,algorithm}

\usepackage{vruler}

\input{commenting.tex}

\allowdisplaybreaks[3]

\graphicspath{ {../} }

\begin{document}

\author[J.~H.~Huggins]{Jonathan H.~Huggins}
\address{Massachusetts Institute of Technology}
\urladdr{http://jhhuggins.org/}
\email{jhuggins@mit.edu}

\author[D.~M.~Roy]{Daniel M.~Roy}
\address{University of Toronto}
\urladdr{http://danroy.org/}
\email{droy@utstat.toronto.edu}

\title[Sequential Monte Carlo as Approximate Samplings]
{Sequential Monte Carlo as Approximate Sampling:
bounds, adaptive resampling via $\infty$-ESS, \\
and an application to Particle Gibbs}

\date{\today}

\begin{abstract}
\input{abstract}

\end{abstract}

\maketitle

\begin{small}
\renewcommand\contentsname{\!\!\!\!}
\setcounter{tocdepth}{1}
\vspace{-15mm}
\tableofcontents
\end{small}

\input{contents}

\bibliographystyle{plainnat}
\bibliography{../main}

\end{document}

%% file: shortcuts2.tex
\usepackage{amsmath}
\usepackage{amsthm}
\usepackage{amsfonts}
\usepackage{amssymb}
\usepackage{dsfont}

\def\[#1\]{\begin{align}#1\end{align}}
\def\(#1\){\begin{align*}#1\end{align*}}

\renewcommand{\var}{\mathbb V}

\newcommand{\as}{\textrm{a.s.}\xspace}

\renewcommand{\Pr}{\mathbb{P}}
\newcommand{\EE}{\mathbb{E}}
\newcommand{\dee}{\mathrm{d}}
\newcommand{\theset}[1]{\lbrace #1 \rbrace}

\newcommand{\kl}[2]{d_{KL}(#1 || #2)}
\newcommand{\totalvar}[2]{d_{TV}(#1, #2)}
\newcommand{\chisq}[2]{d_{\chi^{2}}(#1 || #2)}

\newcommand{\ind}{\mathds{1}}

\newcommand{\ess}{\ensuremath{\mathrm{ESS}}}

\def\opper{operator perspective\xspace}

\def\ciSMC{c$^{i}\alpha$SMC\xspace}
\def\cSMC{c$\alpha$SMC\xspace}
\def\ciARPF{c$^{i}$ARPF\xspace}
\def\cARPF{cARPF\xspace}

\newcommand{\dist}{\sim}

\newcommand{\tuple}[1][]{\left\langle #1 \right\rangle}

\newcommand{\uniformMatrix}{\boldsymbol 1_{1/N}}
\newcommand{\identityMatrix}{\boldsymbol I_{N}}

\newcommand{\truedVanilla}{\pi}
\newcommand{\trued}[1][]{\truedVanilla_{#1}}
\newcommand{\truedRange}[2][1]{\trued[#1:#2]}
\newcommand{\truedAll}[1][t]{\truedRange[1]{#1}}
\newcommand{\truedMarginal}[1][t]{\truedVanilla_{#1}}

\newcommand{\truedPredVanilla}{\eta}
\newcommand{\truedPred}[1][]{\truedPredVanilla_{#1}}
\newcommand{\truedPredRange}[2][1]{\truedPred[#1:#2]}
\newcommand{\truedPredAll}[1][t]{\truedPredRange[1]{#1}}
\newcommand{\truedPredMarginal}[1][t]{\truedPredVanilla_{#1}}

\newcommand{\truedUVanilla}{\gamma}
\newcommand{\truedU}[1][]{\truedUVanilla_{#1}}
\newcommand{\truedURange}[2][1]{\truedU[#1:#2]}
\newcommand{\truedUAll}[1][t]{\truedURange[1]{#1}}

\newcommand{\truedPredUVanilla}{\gamma'}
\newcommand{\truedPredU}[1][]{\truedPredUVanilla_{#1}}
\newcommand{\truedPredURange}[2][1]{\truedPredU[#1:#2]}
\newcommand{\truedPredUAll}[1][t]{\truedPredURange[1]{#1}}

\newcommand{\propd}[1][]{M_{#1}}

\newcommand{\propdAll}[1][t]{\propd[1:#1]}

\newcommand{\propdParam}[2][\theta]{\propd[#2]^{#1}}

\newcommand{\normalizer}[1][t]{Z_{#1}}
\newcommand{\normalizerPred}[1][t]{Z_{#1}'}

\newcommand{\X}[1][]{E_{#1}}
\newcommand{\Xsig}[1][]{\mcE_{#1}}

\newcommand{\Xtsh}{Y}
\newcommand{\Xtshsig}{\mathcal Y}

\newcommand{\range}[3]{#1_{#2:#3}}
\newcommand{\all}[2]{#1_{1:#2}}
\newcommand{\xRange}[2][1]{\range{x}{#1}{#2}}
\newcommand{\yRange}[2][1]{\range{y}{#1}{#2}}
\newcommand{\xAll}[1][t]{\xRange[1]{#1}}
\newcommand{\yAll}[1][t]{\yRange[1]{#1}}

\newcommand{\yiAll}[2][t]{\range{\by}{1}{#1}^{1:#2}}

\newcommand{\potential}[1][t]{g_{#1}}
\newcommand{\potentials}[2][1]{g_{#1:#2}}

\newcommand{\potentialUB}[1][t]{\overline g_{#1}}

\newcommand{\potentialParam}[2][\theta]{\potential[#2]^{#1}}

\makeatletter
\newcommand*\widet[2][0.1ex]{%
        \begingroup
        \mathchoice{\wt@helper{#1}{#2}{\displaystyle}{\textfont}}
                   {\wt@helper{#1}{#2}{\textstyle}{\textfont}}
                   {\wt@helper{#1}{#2}{\scriptstyle}{\scriptfont}}
                   {\wt@helper{#1}{#2}{\scriptscriptstyle}{\scriptscriptfont}}%
        \endgroup
        #2%
}
\newcommand*\wt@helper[4]{%
        \def\currentfont{\the#41}%
        \def\currentskewchar{\char\the\skewchar\currentfont}%
        \setbox\tw@\hbox{\currentfont#2\currentskewchar}%
        \dimen@ii\wd\tw@
        \setbox\tw@\hbox{\currentfont#2{}\currentskewchar}%
        \advance\dimen@ii-\wd\tw@
        \rlap{\raisebox{-#1}{$\m@th#3\kern\dimen@ii\widetilde{\phantom{#2}}$}}%
}
\makeatother

\newcommand{\Samp}[2][n]{X_{#2}^{#1}}
\newcommand{\SampAll}[2][n]{X_{1:#2}^{#1}}
\newcommand{\SampAllPath}[2][n]{\tX_{1:#2}^{#1}}
\newcommand{\SampsT}[1][t]{X_{#1}^{1:N}}
\newcommand{\Samps}[1][t]{\bX_{1:#1}^{1:N}}
\newcommand{\SampsPath}[1][t]{\vphantom{{\tilde X}_{1:#1}^{1:N}}{\smash{\tilde\bX}}_{1:#1}^{1:N}}
\newcommand{\SampAllPathFinal}{\SampAllPath[*]{t}}

\newcommand{\samp}[2][n]{x_{#2}^{#1}}
\newcommand{\sampAll}[2][n]{x_{1:#2}^{#1}}
\newcommand{\sampAllPath}[2][n]{\tx_{1:#2}^{#1}}
\newcommand{\sampsT}[1][t]{x_{#1}^{1:N}}
\newcommand{\samps}[1][t]{\bx_{1:#1}^{1:N}}

\newcommand{\wt}[2][n]{w_{#2}^{#1}}
\newcommand{\wts}[1][t]{w_{#1}^{1:N}}

\newcommand{\Wt}[2][n]{W_{#2}^{#1}}
\newcommand{\Wts}[1][t]{W_{#1}^{1:N}}

\newcommand{\normalizerEst}[1][t]{\hZ_{#1}}
\newcommand{\normalizerEstPred}[1][t]{\hZ'_{#1}}

\newcommand{\truedEst}[1][i]{\truedVanilla^{#1,N}}
\newcommand{\truedEstPred}[1][i]{\truedPredVanilla^{#1,N}}

\newcommand{\truedEstT}[2][i]{\truedEst[#1]_{1:#2}}

\newcommand{\truedEstPredT}[2][i]{\truedEstPred[#1]_{1:#2}}
\newcommand{\truedEstTMarginal}[2][i]{\truedEst[#1]_{#2}}

\newcommand{\truedEstPredTMarginal}[2][i]{\truedEstPred[#1]_{#2}}

\newcommand{\ciEE}[2][i]{\EE^{#1,N}_{#2}}
\newcommand{\ciEEy}[1][i]{\ciEE[#1]{\yiAll{#1}}}
\newcommand{\cEEy}{\ciEE[1]{\yAll}}
\newcommand{\ccEE}[1][\theta]{\ciEE[2]{#1,y^{1:2}}}
\newcommand{\ciKernel}[2][i]{P_{#2}^{#1,N}}
\newcommand{\ciKernely}[1][i]{\ciKernel[#1]{\yiAll{#1}}}
\newcommand{\cKernel}[1][\theta]{\ciKernel[1]{#1,y}}

\newcommand{\ciPr}[2][i]{\Pr_{#2}^{#1,N}}
\newcommand{\ciPry}[1][i]{\ciPr[#1]{\yiAll{#1}}}
\newcommand{\cPry}{\ciPr[1]{\yAll}}

\newcommand{\ancestor}[2][n]{a_{#2}^{#1}}
\newcommand{\ancestorsT}[1][t]{a_{#1}^{1:N}}
\newcommand{\ancestors}[1][t]{\ba_{1:#1}^{1:N}}
\newcommand{\ancestorFinal}[1][t]{\ancestor[*]{#1}}
\newcommand{\Ancestor}[2][n]{A_{#2}^{#1}}
\newcommand{\AncestorsT}[1][t]{A_{#1}^{1:N}}
\newcommand{\Ancestors}[1][t]{\bA_{1:#1}^{1:N}}
\newcommand{\AncestorFinal}{\Ancestor[*]{t}}

\newcommand{\lineage}[2][j]{f_{#2}^{#1}}
\newcommand{\lineagesT}[1][t]{f_{#1}^{1:i}}
\newcommand{\lineages}[1][t]{\bbf_{1:#1}^{1:i}}
\newcommand{\lineageFull}[2][i]{f_{1:#2}^{#1}}
\newcommand{\Lineage}[2][j]{F_{#2}^{#1}}
\newcommand{\LineagesT}[1][t]{F_{#1}^{1:i}}
\newcommand{\Lineages}[1][t]{\bF_{1:#1}^{1:i}}
\newcommand{\LineageFull}[2][i]{F_{1:#2}^{#1}}

\newcommand{\ciNormConst}[2][s]{\mcC_{#1}^{#2}}

\newcommand{\sisE}{\bar\truedVanilla^{\mathrm{S},N}}

\newcommand{\sisET}[1][t]{\sisE_{1,#1}}

\newcommand{\smcCondKernelDist}[1][t]{\range{\tilde\truedVanilla}{1}{#1}}

%% file: commenting.tex
\newcommand{\LATER}[1]{}
\newcommand{\fLATER}[1]{}
\newcommand{\TBD}[1]{}
\newcommand{\fTBD}[1]{}
\newcommand{\PROBLEM}[1]{}
\newcommand{\fPROBLEM}[1]{}

%% file: abstract.tex
Sequential Monte Carlo (SMC) algorithms 
were originally designed for estimating intractable conditional expectations within state-space models,
but are now routinely used to generate approximate samples in the context of general-purpose Bayesian inference.
In particular, SMC algorithms are often used as subroutines within
larger Monte Carlo schemes, and in this context, the demands placed on SMC are different:
control of mean-squared error is insufficient---one needs to control the divergence from the target distribution directly.
Towards this goal, we introduce the \emph{conditional} adaptive resampling particle filter,
building on the work of 
Gordon, Salmond, and Smith~(1993), Andrieu, Doucet, and Holenstein~(2010), and Whiteley, Lee, and Heine~(2016).
By controlling a novel notion of effective sample size, the $\infty$-ESS,  
we establish the efficiency of the resulting SMC sampling algorithm, providing an adaptive resampling extension of the work of
Andrieu, Lee, and Vihola (2013). 
We apply our results to arrive at new divergence bounds for SMC samplers with adaptive resampling
as well as an adaptive resampling version of the Particle Gibbs algorithm with the same 
geometric-ergodicity guarantees as its nonadaptive counterpart.

%% file: contents.tex
\section{Introduction}
\label{sec:introduction}

Sequential Monte Carlo (SMC) methods are a popular class 
of algorithms for approximate inference
\citep{Doucet:2000,Doucet:2001,DelMoral:2006,Kantas:2009,Doucet:2010,Kunsch:2013}.
In the context of Bayesian inference, SMC produces a particle 
approximation to the posterior distribution as well as an unbiased 
estimate of the marginal likelihood.  
Traditionally, particle approximations were built to 
estimate conditional expectations, 
and the analysis of SMC methods focused on this 
\emph{operator} perspective, by bounding the mean squared error of the resulting estimates.

Increasingly, SMC methods are being used to produce approximate samples, 
usually in the inner loop of other approximate inference algorithm.
A key example is the class of particle Markov chain Monte Carlo
(PMCMC) methods, which aim to combine the best features of SMC and MCMC approaches 
by using SMC as a proposal mechanism for a Metropolis--Hastings (``particle MH'')
or approximate Gibbs (``particle Gibbs'') sampler \citep{Holenstein:2009,Andrieu:2010}. 
Characterizing the efficiency of PMCMC methods is an
active area of investigation 
\citep{Andrieu:2009,Andrieu:2014,Andrieu:2013,Chopin:2013b,Lindsten:2014b,Lee:2014}.

When SMC methods are employed for sampling, convergence guarantees
from the \opper are insufficient.
In this work, we take up the \emph{sampling} perspective,
and study the distribution of a sample drawn from the SMC particle filter approximation.
Building off the work of
\citet{Gordon:1993,Andrieu:2010,Andrieu:2013}, and \citet{Whiteley:2013},
we use conditional filters to derive a minorization condition, lower bounding the density of the approximate sample's distribution
in terms of that of the target distribution.  The analysis extends to conditional SMC as well as to adaptive resampling versions.
One of our key contributions is a novel notion of effective sample size, the $\infty$-ESS,  
which we use to establish the efficiency of the adaptive SMC sampling algorithm. 
Thus, our results are both a sampling analogue to the operator work of \citet{Whiteley:2013} 
and an adaptive resampling extension to the sampling work of \citet{Andrieu:2013}.
We apply our results to arrive at new divergence bounds for SMC samplers with adaptive resampling
as well as an adaptive resampling version of the Particle Gibbs algorithm with the same 
geometric-ergodicity guarantees as its nonadaptive counterpart.

In the remainder of this section we provide an overview of our contributions for
the special case of the the conditional adaptive resampling particle filter:
we introduce the conditional adaptive resampling particle filter,
present our main theoretical results characterizing its performance,
and describe an application to a novel adaptive resampling Particle Gibbs algorithm.

\subsection{(Conditional) adaptive resampling particle filters}

We follow a similar setup and notation to \citet{DelMoral:2004}.
Let $(\xi_{t})_{t \ge 1}$ be an inhomogeneous Markov chain on
the measurable space $(\X, \Xsig)$ with transition kernels $(\propd[t])_{t \ge 2}$ and 
initial distribution $\propd[1]$. 
Denote expectations with respect to the Markov chain by $\EE[\cdot]$. 
Let $\potential[t] : \X \to \reals_{+}$, for $t \ge 1$,
be a sequence of $\Xsig$-measurable potential 
functions on $\X$,
let $g_{0} \equiv 1$, and write
$\potentials[s]{t}(\xRange[s]{t}) \defined \prod_{\tau=s}^{t} \potential[\tau](x_{\tau})$.
For $t=1,2,\dots$, define
the measure $\truedAll[t]$ on $\X^t$ given by
\[
\truedAll[t](\dee \xAll[t])
&\defined \truedUAll(\dee \xAll[t]) / \normalizer[t],
\]
where
\[
\truedUAll(\dee \xAll[t]) 
\defined \prod_{s=1}^{t} \potential[s](x_{s})\propd[s](x_{s-1}, \dee x_{s}) 
\quad \text{and} \quad
\normalizer[t] \defined \truedUAll(1). 
\]
(We have written $\propd[1](x_{0},\dee x_1)$ for $\propd[1](\dee x_1)$.) 
Equivalently,
\(
\truedAll[t](\phi) 
&\defined \frac{\EE\left[\phi(\all{\xi}{t})\potentials[1]{t}(\all{\xi}{t})\right]}{\normalizer[t]}, \qquad \phi \colon \X^t \to \reals \text{ measurable},
\)
where $\normalizer[t] \defined \EE\left[\potentials[1]{t}(\all{\xi}{t})\right]$ is the 
normalization constant.\footnote{In the 
state-space setting, the potential $\potential$ would be the 
conditional density (i.e., likelihood) of the observation $v_{t}$ at time $t$ as a function of unobserved state $x_{t}$: i.e.,
$\potential[t](x_{t}) = p_{t}(v_{t} \given x_{t})$. Then $\truedAll[t]$ would be the posterior distribution of the unobserved state sequence
given the observed sequence.}

Towards the goal of efficiently approximating $\truedAll[t]$, we introduce a novel sequential Monte Carlo algorithm:
the $i$-times conditional adaptive resampling particle filter (\ciARPF),
which is a generalization of the adaptive resampling particle filter~\citep{Gordon:1993,Kunsch:2013}
and the conditional SMC algorithm used in particle Gibbs~\citep{Andrieu:2010}. 
(In \cref{sec:condasmc},
we will introduce a further generalization.)
The integer parameter $i \ge 0$ determines
the number of fixed trajectories $\yAll^{1},\dots,\yAll^{i} \in \X^{t}$ required by
the algorithm, which operates by generating a collection $\SampsPath \defined \{ \SampAllPath{t} \}_{n=1}^{N}$
of $N > i$ particles 
with corresponding nonnegative weights $\Wts \defined \{\Wt{t}\}_{n=1}^{N}$.
When $i=0$, we recover the standard (unconditional) adaptive resampling particle filter;
when $i=1$, we recover a generalization of the conditional SMC algorithm that includes adaptive resampling.
For time $s = 1,\dots,t$, the measure $\truedAll[s]$ is approximated by
\(
\truedEstT{s} \defined \sum_{n=1}^{N} \frac{\Wt{s}\potential[s](\Samp{s})}{\sum_{k=1}^{N}\Wt[k]{s}\potential[s](\Samp[k]{s})} \delta_{\SampAllPath[n]{s}} .
\)

The \ciARPF algorithm iteratively constructs $\SampsPath$ and $\Wts$ as follows: 
The first $i$ particles are deterministically set to match the fixed trajectories:
\(
\Samp{s} &= y_{s}^{n}, & \SampAllPath{s} &= \yAll[s]^{n}, &  s=1,\dots,t \text{ and } n=1,\dots,i.
\)
At time $s=1$, the remaining $N-i$ particles $\Samp{1}$, for $n=i+1,\dots,N$, 
are sampled independently and identically from $\propd[1]$.
The corresponding (length $1$) trajectories are
\(
\SampAllPath{1} &= \Samp{1}, & n=i+1,\dots,N.
\)
Furthermore, for all $n=1,\dots,N$, $\Wt{1} = 1$. 

The remaining particle trajectories are generated as follows:
First, we introduce a cutoff parameter $\eta \in [0,1]$ 
and an effective sample size (ESS) function $\ess : \reals_{+}^{N} \to [1,N]$.
The ESS function measures how uniform the current weights $\Wts[s]$ are.
Typically $\ess(\Wts[s]) = 1$ indicates that all but one weight is zero and
$\ess(\Wts[s]) = N$ indicates all the weights are equal. 

For each time $s = 2,\dots,t$:
\bitems

\item If $\ess(\Wts[s-1]) \le \eta N$, a \emph{resampling} step is introduced.
For $n=1,\dots,N$, the weights are set to a common value
\(
\Wt{s} = \Wt[]{s} &\defined \frac 1 N \sum_{k=1}^{N}\Wt[k]{s-1}\potential[s-1](\Samp[k]{s-1})
\)
and, for $n=i+1,\dots,N$, particle $n$'s ``ancestor'' at time $s$, denoted $\Ancestor{s}$, 
is sampled independently, such that $\Ancestor{s} = k$, for $k = 1,\dots,N$, with probability
\(
\frac{\Wt[k]{s-1}\potential[s-1](\Samp[k]{s-1})}{N \,\Wt[]{s}}.
\)

\item If $\ess(\Wts[s-1]) > \eta N$, then 
the algorithm does not resample the particles. For $n=1,\dots,N$, the weights are copied, i.e.,
\(
\Wt{s} &= \Wt{s-1}\potential[s-1](\Samp{s-1}),
\)
and, for $n=i+1,\dots,N$, a record is made that particle $n$ was its own ancestor by setting
 $\Ancestor{s} = n$. 
\item Having sampled the ancestors, the algorithm propagates the particles forward. 
For $n=i+1,\dots,N$,  $\Samp{s}$ is sampled from $\propd[s](\Samp[\Ancestor{s}]{s-1}, \cdot)$,
and the corresponding trajectories are set to
\(
\SampAllPath{s} &= \tuple[\smash{{\SampAll[\Ancestor{s}]{s-1},\Samp{s}}}].
\)
\eitems
In the final step of the algorithm, a single particle $\SampAllPathFinal$
is sampled from the full approximation $\truedEstT{t}$,
and the algorithm yields an estimate of the normalization constant $\normalizer$, 
\(
\normalizerEst \defined \frac{1}{N}\sum_{n=1}^{N} \Wt{t}\potential[t](\Samp{t}).
\)

Let ${\ciEEy}[\cdot]$ denote the expectation operator with
respect to the \ciARPF,
and write
\[\nonumber
\ciKernel{}(\yiAll{i}, \dee \xAll) \defined \ciKernely(\dee \xAll) \defined {\ciEEy}[\delta_{\SampAllPathFinal}(\dee \xAll)]
\]
for the law of $\SampAllPathFinal$ when the $i$ fixed trajectories are $\yiAll{i} \in (\X^{t})^{i}$.
We can now describe in more precise terms how the  \ciARPF kernel $\ciKernel{}$ generalizes several well-known SMC kernels.
When $i=0$, $\truedEstT[0]{t}$ is the standard adaptive SMC particle approximation of $\truedAll[t]$
and $\SampAllPathFinal$ is a single sample from the particle approximation.
When $i=1$ and resampling is done at every step by taking $\eta=1$,
$\ciKernel[1]{}$ %
is exactly the conditional SMC kernel used 
in particle Gibbs samplers~\citep{Holenstein:2009,Andrieu:2010,Andrieu:2013}.
For general $\eta \in (0,1)$, we obtain a novel adaptive resampling variant that we study in the sequel.
In particular, under mild regularity conditions, 
$\ciKernel[1]{}$ %
defines a Markov kernel with invariant distribution $\truedAll[t]$.

\subsection{Controlling c$^{i}$ARPF efficiency with $\infty$-ESS}

We can analyze the quality of the c$^{i}$ARPF kernel $\ciKernel{}$ by 
quantifying the extent to which high-probability sets under the target 
distribution also have high probability under the kernel. 
The following theorem 
establishes a minorization condition  
for the $i$-times conditional filter in terms of the $(i+1)$-times conditional filter:

\bthm \label{thm:ci-arpf-rd-general}
For all $t \ge 1$, $i \ge 0$, $N > i$, and $\yAll^{1},\dots,\yAll^{i} \in \X^{t}$, 
\[
\ciKernel{}(\yiAll{i}, S)
&\ge (1-i/N)^t \int_{S}\frac{\normalizer}{\ciEEy[i+1][\normalizerEst]} \truedAll[t](\dee \yAll^{i+1}), 
& S \subseteq \X^{t} \text{ measurable}. \label{eq:main-result}
\]
\ethm

The integral appearing in \cref{eq:main-result} has no simple form in general, but in many 
settings we will be able to obtain a lower bound on the integrand that does not depend on the 
fixed trajectories $\yiAll{i+1}$.
In those cases, the integral is simply replaced by this uniform lower bound.
For $i=0$, we are then immediately able to control numerous measures of divergence between $\truedAll[t]$ and 
$\ciKernel[0]{}$ (i.e., the law of $\SampAllPathFinal$). 
For example, in the case of total variation distance, we have the following corollary to \cref{thm:ci-arpf-rd-general}: 
\bcor \label{cor:arpf-kl}
If $\ciEE[1]{\yAll}[\normalizerEst/\normalizer] \le B_{t,N}$ for all $\yAll \in \X^{t}$, then
\[
\totalvar{\truedAll[t]}{\ciKernel[0]{}} \le 1 - B_{t,N}^{-1}.
\]
\ecor
For $i=1$, a uniform lower bound assumption implies a minorization condition on
the kernel $\ciKernel[1]{}(\yAll, \dee \xAll)$, which in turn implies fast
mixing of the Markov chain with kernel $\ciKernel[1]{}$:
\bcor \label{cor:carpf-tv}
If $\ciEEy[2][\normalizerEst/\normalizer] \le B_{t,N}$ for all $\yAll^{1}, \yAll^{2} \in \X^{t}$, then
the Markov chain with transition kernel $\ciKernel[1]{}(\yAll, \dee \xAll)$ is uniformly ergodic in
total variation distance and has invariant distribution $\truedAll$.
In particular, for all $\yAll \in \X^{t}$ and $k \ge 1$,
\[
\totalvar{\tilde\pi_{\yAll}^{k}}{\truedAll} &\le \left(1 - B_{t,N}^{-1}(1-1/N)^{t}\right)^{k},
\]
where $\tilde\pi_{\yAll}^{k} \defined \delta_{\yAll}[\ciKernel[1]{}]^{k}$ is the law of the Markov chain, with initial state $\yAll$, after $k$ transitions.
\ecor

In order to apply the corollaries, it remains to bound
${\ciEEy}[\normalizerEst/\normalizer]$.
Such a bound was obtained for the \emph{nonadaptive} conditional SMC kernel in \citet{Andrieu:2013}.
However, in our general adaptive resampling setting, one must make a careful choice of 
effective sample size function. 
To this end, we introduce a generalized notion of effective sample size,
which includes several existing definitions as special cases.
For $p \in (1,\infty]$, let $p_{*} \defined \frac{p}{p - 1}$ be the
conjugate exponent of $p$ (so $1/p + 1/p_{*} = 1$). 
The \textbf{$p$-effective sample size} ($p$-ESS) of the weight vector 
$w^{1:N} \in \reals_{+}^{N}$ is
\(
\ess_{p}(w^{1:N})
&\defined \frac{\|w^{1:N}\|_{1}^{p_{*}}}{\|w^{1:N}\|_{p}^{p_{*}}}.
\)
The following proposition highlights some elementary properties of $p$-ESS.

\bprop \label{prop:basic-ess-properties}
The $p$-ESS has the following properties:
\benum[leftmargin=*]
\item For all $p \in (1,\infty]$, $1 \le \ess_{p}(w^{1:N}) \le N$. 
The lower bound is achieved if and only if all but one of the weights is zero. 
The upper bound is achieved if and only if all the weights are equal.
\item For $1 < p < q \le \infty$, $\ess_{p}(w^{1:N}) \ge \ess_{q}(w^{1:N}) \ge N^{-(1 - q_{*}/p_{*})}\ess_{p}(w^{1:N})$,
with equality if and only if $K$ weights are equal and the rest are zero.
\eenum
\eprop

Part (1) demonstrates that the $p$-ESS satisfies basic properties one would expect of a measure of effective sample size.
Part (2) places the family of $p$-ESS measures in a linear order: the larger the value of $p$, the more stringent the notion of effective
sample size. 

The standard definition of effective sample size is precisely the $2$-ESS.
\citet{Whiteley:2013} provided a rigorous 
justification for the use of 2-ESS from the \opper:
if adaptive resampling is used to guarantee that the 2-ESS does not fall below $\zeta N$, 
for some fixed parameter $\zeta \in (0,1]$, 
then the error bounds on the operator approximation match those of the nonadaptive sampler with $\zeta N$ particles. 
More formally, let $\Im_{t}(\yAll) = y_{t}$ be the projection onto the $t$-th component.
Under appropriate regularity conditions,
for every bounded measurable $\phi : \X \to \reals$ and real $r \ge 1$,
\(
\sup_{t \ge 1} \ess_{2}(\Wts) \ge \zeta N
\implies 
\sup_{t \ge 1} \ciEE[0]{}\left[|\truedEstT[0]{t}(\phi \circ \Im_{t}) - \trued(\phi \circ \Im_{t})|^{r}\right]^{1/r} 
	\le \frac{a(r)b(\phi)}{\sqrt{\zeta N}},
\)
where $a(r)$ and $b(\phi)$ are explicit functions. 

To upper bound ${\ciEEy}[\normalizerEst/\normalizer]$, however, we will require a lower
bound on the $\infty$-ESS, which by \cref{prop:basic-ess-properties}(2) is a more stringent
notion of effective sample size than 2-ESS:

\begin{restatable}{assumption}{infessbound} \label{asm:inf-ess-bound}
There exists $\zeta \in (0,1]$ such that $\ess_{\infty}(\Wts[s]) \ge \zeta N$
for all $1 \le s \le t$. 
\end{restatable}

For the choice
$\ess = \ess_{\infty}$ (i.e., under \cref{asm:inf-ess-bound}),
we can bound the estimate of the normalization constant.
Let $G_{s,t}(x_{s}) \defined \EE[\potentials[s]{t-1}(\range{\xi}{s}{t-1}) \given \xi_{s} = x_{s}]$
for $s=1,\dots,t$ and $G_{0,t} \defined \EE[\potentials[1]{t-1}(\all{\xi}{t-1})]$.
We now arrive at our second main result:

\bthm \label{thm:ciEE-bound}
If \cref{asm:inf-ess-bound} holds, then for all $i,t \ge 1$, $N > i$, 
$\yAll^{1},\dots,\yAll^{i} \in \X^{t}$, 
\[
{\ciEEy}[\normalizerEst/\normalizer]
&\le 1 + \frac{\normalizer^{-1}\sum_{s=1}^{t}\sum_{j=1}^{i}G_{0,s}G_{s,t+1}(\yAll^{i}) - \zeta i}{\zeta N}  + \Theta(N^{-2}).
\]
\ethm

Two possible further assumptions both lead to uniform bounds on ${\ciEEy}[\normalizerEst]$.

\begin{restatable}{assumption}{boundedpotentials} \label{asm:bounded-potentials}
The potentials satisfy $\potentialUB[s] \defined \sup_{x \in \X} \potential[s](x) < \infty$ 
for all $1 \le s \le t$.
\end{restatable}
\begin{restatable}{assumption}{Gzeroratiobound} \label{asm:G0-ratio-bound}
There exists a constant $\beta > 0$ such that for any $t, s \in \naturals$,
\(
\sup_{x \in \X} \frac{G_{0,t}G_{t,t+s}(x)}{G_{0,t+s}} \le \beta.
\)
\end{restatable}

\bcor \label{cor:ci-arpf-quantitative-bounds}
Under the same conditions as \cref{thm:ciEE-bound},
if \cref{asm:bounded-potentials} holds then
\[
{\ciEEy}[\normalizerEst/\normalizer] 
\le 1 + \normalizer[t]^{-1}\prod_{s=1}^{t}\potentialUB[s]\left[\left(1 + \frac{i}{\zeta N}\right)^{t} - 1\right]
\label{eq:normalizer-simple-bound}
\]
while if \cref{asm:G0-ratio-bound} holds then 
\[
{\ciEEy}[\normalizerEst/\normalizer] 
\le  \left(1 + \frac{\beta}{\zeta N}\right)^{t}.
\label{eq:normalizer-beta-bound}
\]
\ecor

Combining \cref{cor:arpf-kl,cor:ci-arpf-quantitative-bounds} yields the following 
guarantees for the ARPF sampler:
\bthm
If \cref{asm:inf-ess-bound,asm:bounded-potentials} hold then
\[
\totalvar{\truedAll[t]}{\ciKernel[0]{}} 
\le \frac{t\normalizer[t]^{-1}\prod_{s=1}^{t}\potentialUB[s]}{\zeta N + t\normalizer[t]^{-1}\prod_{s=1}^{t}\potentialUB[s]} + \Theta(N^{-2})
\]
while if \cref{asm:inf-ess-bound,asm:G0-ratio-bound} hold then
\[
\totalvar{\truedAll[t]}{\ciKernel[0]{}} 
\le  \frac{\beta t}{\zeta N + \beta t} + \Theta(N^{-2}).
\]
\ethm

\subsection{Applications to Particle Gibbs}
\label{sec:intro-applications}

In the language of state-space models, the setting described so far involves 
approximating the posterior distribution of a Markov chain given indirect stochastic 
observations of the chain's values. 
However, it is often the case that the chain and the potentials are controlled by a global
parameter $\theta \in \Theta$ for which there is a prior distribution $\varpi(\dee\theta)$. 
Replace $\propd[s]$ by $\propdParam{s}$ and 
$\potential[s]$ by $\potentialParam{s}$, then parameterize
the other quantities defined previously in terms of $\propd[s]$ 
and $\potential[s]$ by $\theta$. 
Let $(\Xtsh, \Xtshsig) \defined (\X^{t}, \mcB(\X^{t}))$.
We will suppress much of the time dependence when possible to make the
notation less cluttered.
The target distribution on the product space 
$(\Theta \times \Xtsh, \mcB(\Theta \times \Xtsh))$ is
\[
\trued(\dee \theta \times \dee y)
&\defined \truedU(\dee \theta \times \dee y) / \normalizer[],
\]
where
\[
\truedU(\dee \theta \times \dee y) 
\defined \varpi(\dee \theta)\prod_{s=1}^{t} \potentialParam{s}(y_{s})\propdParam{s}(y_{s-1}, \dee y_{s}) 
\quad \text{and} \quad
\normalizer[] \defined \truedU(1). 
\]
Let $\trued[\theta](\dee y)$ and $\trued[y](\dee \theta)$
denote the disintegrations of $\pi$ along $\Theta$ and along $Y$, respectively.

The particle Gibbs sampler approximates the two-stage Gibbs kernel
\[
\Pi(\theta, y, \dee \vartheta \times \dee z)
\defined \trued[y](\dee \vartheta)\trued[\vartheta](\dee z).
\]
In many settings, such as non-linear or non-Gaussian
state-space models, it is possible to sample 
from $\trued[y](\dee \vartheta)$, but difficult to
sample from $\trued[\vartheta](\dee z)$.
The idea is to replace $\trued[\vartheta](\dee z)$
with an SMC-based approximation $\Pi_{\vartheta}(y, \dee z)$
that leaves $\trued[\vartheta](\dee z)$ invariant,
leading to a kernel of the form
$\trued[y](\dee \vartheta)\Pi_{\vartheta}(y, \dee z)$.

We introduce the adaptive resampling particle Gibbs (ARPG) sampler, 
which employs the \cARPF kernel $\cKernel[\vartheta](\dee z)$ 
to approximate the conditional distribution $\trued[\vartheta](\dee z)$ 
that would be used in a standard Gibbs sampler. 
The ARPG kernel is thus given by
\[
\Pi_{N}(\theta, y, \dee \vartheta \times \dee z)
\defined \trued[y](\dee \vartheta)\cKernel[\vartheta](\dee z).
\]
\cref{thm:ci-arpf-rd-general,cor:carpf-tv,cor:ci-arpf-quantitative-bounds}
together with the results of \citet{Andrieu:2013}
yield guarantees on the ergodicity properties of the \cARPF kernel and the ARPG sampler.
Once instances of $N$ are replaced by $\zeta N$, the guarantees essentially
match those provided by \citet{Andrieu:2013} for the standard PG sampler.

\bthm
If \cref{asm:inf-ess-bound} holds, then the \cARPF kernel and ARPG sampler have 
the following properties:
\benum[leftmargin=*]
\item If \cref{asm:bounded-potentials} holds then there exists
$\eps_{t,N} = 1 - C_{t}/N$ such that for any $\theta \in \Theta$, $y \in \Xtsh$, and $k \ge 1$, 
\[
d_{TV}(\delta_{y}[\ciKernel[1]{\theta}]^{k}, \trued[\theta]) \le (1 - \eps_{t,N})^{k}. \label{eq:carpf-uniform-ergodicity}
\]
\item If \cref{asm:G0-ratio-bound} holds and $N \ge t/C + 1$ for any fixed $C > 0$, 
then for any $t \ge 1$, \cref{eq:carpf-uniform-ergodicity} holds with
\[
\eps_{t,N} \ge \exp\left(-\frac{C}{\zeta}(2\beta + \zeta)\right).
\]
\item If either \cref{asm:bounded-potentials} or \cref{asm:G0-ratio-bound} holds,
then whenever the Gibbs sampler is geometrically ergodic the ARPG sampler is 
geometrically ergodic as well. 
\eenum
\ethm

At a high level, the results we have obtained highlight the role of the
expected value of $\normalizerEst$ 
in the mixing properties of conditional SMC Markov chains and particle Gibbs (PG) samplers:
In order to show geometric ergodicity for adaptive resampling particle Gibbs samplers, 
it suffices to establish bounds on the expected value of $\normalizerEst$ under the twice-conditional filter, 
and the growth of the expectation as $t$ increases determining how well the particle Gibbs 
algorithm scales.
Similarly, a bound on the expected value of $\normalizerEst$ under the once-conditional filter
implies a bound on $\totalvar{\truedAll}{\ciKernel[0]{}}$.
Hence, as a slogan, good performance of (adaptive resampling) particle Gibbs is equivalent 
to good performance of (adaptive) SMC for sampling.

\section{Preliminaries}
\label{sec:smc}

In this section, we fix some additional notation, 
introduce a few key additional definitions,
and then present $\alpha$SMC~\citep{Whiteley:2013}, 
a generalization of the adaptive resampling particle filter described in
the introduction. 

For a positive integer $K$, let $[K] \defined \theset{1,2,\dots,K}$. 
If $x_{i},\dots,x_{j}$ are elements of a sequence, write 
$\xRange[i]{j} \defined \tuple[x_{i},x_{i+1},\dots,x_{j}]$.
We use the following conventions: $\sum_{\emptyset} = 0$,
$\prod_{\emptyset} = 1$, and $0/0 = 0$. 

Let $(S, \mcS)$, $(S', \mcS')$ be measurable spaces. 
Then $K : S \times \mcS' \to \reals$ is a kernel if 
$K(\cdot,B)$ is a $(S,\mcS)$-measurable function for 
all $B \in \mcS'$ and $K(x, \cdot)$ is measure on $(S', \mcS')$
for all $x \in \mcS$. 
For a measure $\mu$ on $(S, \mcS)$
and kernels $K, K' : S \times \mcS \to \reals$,
let $\mu K(\dee y) \defined \int \mu(\dee x)K(x, \dee y)$ 
and $KK'(x, \dee z) \defined \int K(x, \dee y)K'(y, \dee z)$. 
We will often use measures and kernels as operators.
For a measurable function $\phi : S \to \reals$, 
let $\mu(\phi) \defined \EE_{\xi \dist \mu}[\phi(\xi)] = \int \phi(x)\mu(\dee x)$ and
$K(\phi)(x) \defined \int \phi(y)K(x, \dee y)$.
For measures $\mu, \nu$ on $(S,\mcS)$, we will write
$\mu \ll \nu$ to denote that $\mu$ is absolutely continuous 
with respect to $\nu$, in which case we will write 
$\dee \mu/\dee \nu$ for the $\nu$-almost everywhere 
($\nu$-a.e.) unique function $f$ satisfying 
$\mu(A) = \int_A f \,\dee \nu$, for all $A \in \mcS$.
When the choice is clear from context, we may
write $\mcB(S)$ for the $\sigma$-algebra of the 
space $S$. 

For probability measures $\mu$ and $\nu$ on $(S, \mcS)$,
the total variation distance between
$\mu$ and $\nu$ is 
\[
d_{TV}(\mu, \nu) \defined \sup_{A \in \mcS} |\mu(A) - \nu(A)|.
\]
If $\mu \ll \nu$, then the Kullback--Liebler (KL) divergence is 
\[
\kl{\mu}{\nu} \defined \mu(\log \dee\mu/\dee\nu) \label{eq:kl-def}
\]
and the $\chi^{2}$ divergence is
\[
d_{\chi^{2}}(\mu||\nu) 
\defined \nu([\dee\mu/\dee\nu - 1]^{2})
= \mu(\dee\mu/\dee\nu) - 1.
\]

Finally, we note that, when there is little risk of confusion,
we will ignore measure-theoretic
niceties such as the distinction between equality 
and a.e.-equality.

Recall that $(\xi_{t})_{t \ge 1}$ is an inhomogeneous Markov chain on
$(\X, \Xsig)$ with transition kernels $(\propd[t])_{t \ge 2}$ and 
initial distribution $\propd[1]$,
and that $\EE[\cdot]$ denotes expectation with respect to the Markov chain.
We will write $\propd[1](x_{0},\cdot)$ for $\propd[1](\cdot)$ when convenient
and,
for all $t \ge 1$ and $x_{t-1} \in \X$, we will assume that
$\propd[t](x_{t-1},\cdot)$ has a density with respect to
some common $\sigma$-finite dominating measure (which we 
denote by $\dee x$). 
We will abuse notation and write 
$\propd[t](x_{t-1},x_{t})$
for
the density of 
$\propd[t](x_{t-1},\cdot)$ as $x_{t}$.
Recall that, for each $t \ge 1$, 
$\potential[t] : \X \to \reals_{+}$
denotes a 
$\Xsig$-measurable potential function,
with $g_{0} \equiv 1$. Finally, recall that
$\potentials[s]{t}(\xRange[s]{t}) \defined \prod_{\tau=s}^{t} \potential[\tau](x_{\tau})$.

\subsection{Target distributions}

We now introduce some additional target distributions. 
(We will also repeat the definition of $\truedUAll[t]$ and  $\truedAll[t]$ for completeness.)

Let $\all{\phi}{t} \colon \X^t \to \reals$ and $\phi_{t} \colon \X \to \reals$ denote generic measurable functions. 
For each $t \ge 1$, the unnormalized predictive and updated measures are
defined, respectively, by
\(
\truedPredUAll[t](\all{\phi}{t}) 
\defined \EE\left[\all{\phi}{t}(\all{\xi}{t})\potentials[1]{t-1}(\all{\xi}{t-1})\right]
\qquad \text{and} \qquad 
\truedUAll[t](\all{\phi}{t}) 
\defined \EE\left[\all{\phi}{t}(\all{\xi}{t})\potentials[1]{t}(\all{\xi}{t})\right]
\)
with corresponding marginal versions
\(
\truedPredU[t](\phi_{t}) 
\defined \EE\left[\phi_{t}(\xi_{t})\potentials[1]{t-1}(\all{\xi}{t-1})\right]
\qquad \text{and} \qquad 
\truedU[t](\phi_{t}) 
\defined \EE\left[\phi_{t}(\xi_{t})\potentials[1]{t}(\all{\xi}{t})\right].
\)
Our ultimate goal is to approximate the normalized predictive and updated measures
along with their marginal versions:
\[
\truedAll[t](\all{\phi}{t}) 
&\defined \frac{\truedUAll[t](\all{\phi}{t})}{\normalizer[t]},
&\truedPredAll[t](\all{\phi}{t}) 
&\defined \frac{\truedPredUAll[t](\all{\phi}{t})}{\normalizerPred[t]}, \\
\trued[t](\phi_{t}) &\defined \frac{\truedU[t](\phi_{t})}{\normalizer[t]}, 
\qquad\qquad \text{ and}
& \truedPred[t](\phi_{t}) &\defined \frac{\truedPredU[t](\phi_{t})}{\normalizerPred[t]},
\]
where $\normalizer[t] \defined \truedU[t](1)$ and 
$\normalizerPred[t] \defined \truedPredU[t](1)$ are 
normalization constants.

\subsection{The $\alpha$SMC algorithm}
\label{sec:smc-algorithms}

In the introduction, adaptation in the particle filter 
was implemented via a simple multinomial resampling step, triggered when the effective sample size fell below a fixed threshold.
For the remainder of the article, we will consider a more general mechanism for adaptation 
captured by the $\alpha$SMC algorithm introduced by \citet{Whiteley:2013}.
The $\alpha$SMC algorithm can produce 
sequential importance sampling (SIS), sampling importance 
resampling (SIR, also known as the bootstrap filter), 
and numerous other SMC variants as special cases. 
Not only does the $\alpha$SMC formulation aid in 
analyzing adaptive resampling strategies, it provides
a useful framework for devising novel adaptive schemes
with attractive computational properties, such as
admitting parallelization even on resampling steps. 
In the remainder of this section, we outline the (unconditional) $\alpha$SMC algorithm.
In the following section, we introduce a novel $i$-times conditional version of $\alpha$SMC, which will include the \ciARPF as a special case. 

The $\alpha$SMC algorithm, which is given as \cref{alg:alpha-smc},
provides a flexible resampling mechanism: at each time $t$, a stochastic matrix 
$\alpha_{t-1}$ is chosen from a set $\mathbb A_{N}$ of $N \times N$ matrices. 
We denote the value in the $n$-th row and 
$k$-th column of $\alpha_{t-1}$ by $\alpha_{t-1}^{nk}$. 
The $\alpha$SMC estimators are
\(
\truedEstT[0]{t} \defined \sum_{n=1}^{N} \frac{\Wt{t}\potential[t](\Samp{t})}{\sum_{k=1}^{N}\Wt[k]{t}\potential[t](\Samp[k]{t})} \delta_{\SampAllPath[n]{t}},
\qquad
\truedEstTMarginal[0]{t} \defined \sum_{n=1}^{N} \frac{\Wt{t}\potential[t](\Samp{t})}{\sum_{k=1}^{N}\Wt[k]{t}\potential[t](\Samp[k]{t})} \delta_{\Samp[n]{t}},
\)
\(
\truedEstPredT[0]{t} \defined \sum_{n=1}^{N} \frac{\Wt{t}}{\sum_{k=1}^{N}\Wt[k]{t}} \delta_{\SampAllPath[n]{t}},
\quad
\text{and} \quad 
\truedEstPredTMarginal[0]{t} \defined \sum_{n=1}^{N} \frac{\Wt{t}}{\sum_{k=1}^{N}\Wt[k]{t}} \delta_{\Samp[n]{t}},
\)
and the estimators of the normalization constants $\normalizer$
and $\normalizerPred[t]$ are
\(
\normalizerEst \defined \frac{1}{N}\sum_{n=1}^{N} \Wt{t}\potential[t](\Samp{t})
\quad \text{and} \quad
\normalizerEstPred \defined \frac{1}{N}\sum_{n=1}^{N} \Wt{t}.
\)
Expectations with respect the law of the $\alpha$SMC algorithm
are written as $\ciEE[0]{}[\cdot]$.

SIS, SIR, and the standard adaptive algorithm are
obtained as special cases of $\alpha$SMC as follows.
SIS is recovered when $\alpha_{t-1} = \identityMatrix$, 
the $N \times N$ identity matrix, while
SIR is recovered when $\alpha_{t-1} = \uniformMatrix$, 
the $N \times N$ matrix with all entries equal to $1/N$. 
The adaptive particle filter (APF) algorithm is 
obtained by setting $\alpha_{t-1}$ to $\uniformMatrix$ 
if $\ess_{2}(\Wts[t-1]) \le \zeta N$ and to $\identityMatrix$ otherwise,
where $\zeta \in (0,1]$ is fixed.

\begin{algorithm}[t]
\caption{$\alpha$SMC}
\label{alg:alpha-smc}
\begin{algorithmic}
\FOR{$n = 1,\dots,N$}
	\STATE Sample $\Samp{1}$ from $\propd[1]$
	\STATE Set $\SampAllPath{1} \gets \Samp{1}$
	\STATE Set $\Wt{1} \gets 1$
\ENDFOR
\FOR{$t=2,3,\dots$}
	\STATE Select $\alpha_{t-1}$ from $\mathbb A_{N}$ according to some function of $\Samps[t-1]$ and $\Ancestors[t-2]$
	\FOR{$n = 1,\dots,N$}
			\STATE Set $\Wt{t} \gets \sum_{k=1}^{N} \alpha_{t-1}^{nk}\Wt[k]{t-1}\potential[t-1](\Samp[k]{t-1})$
		\STATE Sample $\Ancestor{t-1}$ from $\distMulti\left(\tuple[{\frac{ \alpha_{t-1}^{nk}\Wt[k]{t-1}\potential[t-1](\Samp[k]{t-1})}{\Wt{t}}}]_{k=1}^{N}\right)$
		\STATE Sample $\Samp{t}$ from $\propd[t](\Samp[\Ancestor{t-1}]{t-1}, \cdot)$
		\STATE Set $\SampAllPath{t} \gets \tuple[{\SampAll[\Ancestor{t-1}]{t-1},\Samp{t}}]$
	\ENDFOR
\ENDFOR
\end{algorithmic}
\end{algorithm}

\begin{algorithm}[t]
\caption{Conditional $\alpha$SMC}
\label{alg:c-alpha-smc}
\begin{algorithmic}
\REQUIRE Fixed trajectory $\yAll$
\STATE Sample $\Lineage[1]{1}$ uniformly from $\{1,\dots,N\}$
\FOR{$n = 1,\dots,N$}
	\IF{$n = \Lineage[1]{1}$}
		\STATE Set $\Samp{1} \gets y_{1}$
	\ELSE
		\STATE Sample $\Samp{1}$ from $\propd[1]$
	\ENDIF
	\STATE Set $\SampAllPath{1} \gets \Samp{1}$
	\STATE Set $\Wt{1} \gets 1$
\ENDFOR
\FOR{$t=2,3,\dots$}
	\STATE Select $\alpha_{t-1}$ from $\mathbb A_{N}$ according to some function of $\Samps[t-1]$  and $\Ancestors[t-2]$
	\STATE Sample $\Lineage[1]{t}$ from $\distMulti\left(\tuple[\alpha^{k\Lineage[1]{t-1}}_{t-1}]_{k=1}^{N}\right)$
	\FOR{$n = 1,\dots,N$}
		\STATE Set $\Wt{t} \gets \sum_{k=1}^{N} \alpha_{t-1}^{nk}\Wt[k]{t-1}\potential[t-1](\Samp[k]{t-1})$
		\IF{$n = \Lineage[1]{t}$}
			\STATE Set $\Ancestor{t-1} \gets \Lineage[1]{t-1}$ 
			\STATE Set $\Samp{t} \gets y_{t}$
		\ELSE
			\STATE Sample $\Ancestor{t-1}$ from $\distMulti\left(\tuple[{\frac{ \alpha_{t-1}^{nk}\Wt[k]{t-1}\potential[t-1](\Samp[k]{t-1})}{\Wt{t}}}]_{k=1}^{N}\right)$
			\STATE Sample $\Samp{t}$ from $\propd[t](\Samp[\Ancestor{t-1}]{t-1}, \cdot)$
		\ENDIF
		\STATE Set $\SampAllPath{t} \gets \tuple[{\SampAll[\Ancestor{t-1}]{t-1},\Samp{t}}]$
	\ENDFOR
\ENDFOR
\STATE Sample $\AncestorFinal$ from $\distMulti\left(\tuple[{\frac{\Wt[k]{t}\potential[t](\Samp[k]{t})}{\sum_{n=1}^{N}\Wt{t}\potential[t](\Samp{t})}}]_{k=1}^{N}\right)$
\end{algorithmic}
\end{algorithm}

\section{Conditional $\alpha$SMC} 
\label{sec:condasmc}

It is useful both algorithmically and analytically to  
generalize $\alpha$SMC in such a way that one or more particle trajectories
is fixed ahead of time.
The result, which we will refer to as
\emph{conditional} $\alpha$SMC,
is a strict generalization of the conditional adaptive particle filter given in the introduction.
We will use conditional $\alpha$SMC 
to study the properties of (unconditional) $\alpha$SMC,
to design novel adaptive particle Gibbs algorithms,
and 
to analyze their mixing properties.

For this section, fix $t \ge 1$, $i \ge 0$, and $N > i$. 
The \emph{$i$-times conditional $\alpha$SMC
(\ciSMC) process} (or simply the \emph{\cSMC process} when $i=1$) is defined on the space 
$(\X^{N} \times [N]^{N} \times [N]^{i})^{t-1} \times \X^{N} \times [N] \times [N]^{i}$,
and is essentially equivalent to $\alpha$SMC with the first $i$ particle trajectories,
but not their lineages,  fixed a priori. 
If $\lineagesT[] \in [N]^{i}$ are indices of the first $i$ particles, let 
$\mcD(\lineagesT[]) \defined \prod_{j \ne j'} \ind(\lineage[j]{} \ne \lineage[j']{})$
be the function that indicates whether the indices
are distinct. 
As in $\alpha$SMC, the matrix $\alpha_{t-1} \in \mathbb A_{N}$ is a function 
of $\Samps[t-1]$ and $\Ancestors[t-2]$.
We have $\samps \in (\X^{N})^{t}$, $\lineages \in ([N]^{i})^{t}$,
$\ancestors[t-1] \in ([N]^{N})^{t-1}$, and $\ancestorFinal \in [N]$, and
use the notation
\[
\wt{1} &\defined 1,  &
\wt{t} &\defined \sum_{k=1}^{N} \alpha_{t-1}^{nk}\wt[k]{t-1}\potential[t-1](\samp[k]{t-1}),
\]
and
\[
r_{n}(k | \wts[s-1], \samps[s-1]) 
\defined \frac{\alpha_{s-1}^{nk}\wt[k]{s-1}\potential[s-1](\samp[k]{s-1})}{\wt{s}}.
\]
For fixed trajectories $\yAll^{1},\dots,\yAll^{i} \in \X^{t}$, 
the law of the \ciSMC process is given by
\[
{\ciPry}[\SampsT[1] \in \dee \sampsT[1], \LineagesT[1] = \lineagesT[1]]
 &\defined \ciNormConst[1]{i} \mcD(\lineagesT[1])\prod_{j=1}^{i}\frac{1}{N}\delta_{y_{1}^{j}}(\dee \samp[{\lineage[j]{1}}]{1})
  \prod_{n \notin \lineagesT[1]}^{N} \propd[1](\dee \samp{1}),
\]
for $s = 2,\dots,t$,
\[
\begin{split}
&{\ciPry}[\SampsT[s] \in \dee \sampsT[s], \AncestorsT[s-1] = \ancestorsT[s-1], \LineagesT[s] = \lineagesT[s] \given \\
&\phantom{\ciPry(} \Samps[s-1] = \samps[s-1], \Ancestors[s-2] = \ancestors[s-2], \LineagesT[s-1] = \lineagesT[s-1]]  \\
  &\defined \ciNormConst{i} \mcD(\lineagesT[s]) \prod_{j=1}^{i}\alpha_{s-1}^{\lineage[j]{s}\lineage[j]{s-1}} \delta_{y_{s}^{j}}(\dee \samp[{\lineage[j]{s}}]{s})\ind(\ancestor[{\lineage[j]{s}}]{s-1} = \lineage[j]{s-1}) \\
  &\phantom{\defined~} \times \prod_{n \notin \lineagesT[s]} r_{n}(\ancestor{s-1} | \wts[s-1], \samps[s-1])\propd[s](\samp[\ancestor{s-1}]{s-1}, \samp{s})
\end{split}
\]
and
\[
{\ciPry}[\AncestorFinal = \ancestorFinal \given \Samps = \samps, \Ancestors[t-1] = \ancestors[t-1]]
&\defined \frac{\wt[\ancestorFinal]{t} \potential(\samp[\ancestorFinal]{t})}{\sum_{n=1}^{N} \wt{t}\potential(\samp{t})}.
\]
The $\ciNormConst{i}$ terms are normalization constants 
that ensure the expressions are valid probabilities.
Let $\SampAllPathFinal \defined \SampAll[\AncestorFinal]{t}$,
let ${\ciEEy}[\cdot]$ denote the expectation operator with
respect to the \ciSMC, and write
\[\nonumber
\ciKernel{}(\yiAll{i}, \dee \xAll) \defined \ciKernely(\dee \xAll) \defined {\ciEEy}[\delta_{\SampAllPathFinal}(\dee \xAll)]
\]
for the law of $\SampAllPathFinal$.

The normalization constants $\ciNormConst{i}$ arise because the lineages $\lineages$ 
of the fixed trajectories $\yiAll{i}$ are kept distinct.
The $c^{i}\alpha$SMC kernel enforces distinct lineages for the fixed trajectories 
since in general $\yAll^{j} \ne \yAll^{j'}$ for $j \ne j'$ and, from an algorithmic 
standpoint, allowing overlapping lineages could lead to a substantial increase in 
complexity, both in terms of implementation and computation. 
The distinct lineage requirement is enforced by the $\mcD(\lineagesT[s])$ terms.
Since there is at most one fixed trajectory when $i=0$ or 1, 
$\ciNormConst{0} = \ciNormConst{1} = 1$ for all $s \in [t]$. 

\cref{alg:c-alpha-smc} provides pseudocode to sample from the law of the \cSMC process, 
which is a necessary part of the particle Gibbs sampler described in \cref{sec:particle-mcmc}.
Sampling from the law of the $c^{i}\alpha$SMC process for $i > 1$ is more delicate, but 
unnecessary since these are only used for analytical purposes. 

\begin{remark}
To recover the \ciARPF described in the introduction, let $\alpha_{t-1} = \uniformMatrix$ 
if $\ess(\Wts[t-1]) \le \zeta N$ and let $\alpha_{t-1} = \identityMatrix$ otherwise.
Then note that by the symmetry of $\uniformMatrix$ and $\identityMatrix$, instead of sampling them,
we can set the lineage for the $j$-th fixed trajectory to $j$: that is, set
$\lineage[j]{s} = j$ for all $j \in [i]$ and $s \in [t]$.
\end{remark}

\section{Main Results}
\label{sec:main}

We are now ready to undertake our study of the $i$-times conditional $\alpha$SMC kernel
$\ciKernel{}(\yiAll{i}, \dee \xAll)$. 
Specifically, our aim is to understand the conditions under which 
the \ciSMC kernel is close to $\truedAll(\dee \xAll)$. 
Formally, we will establish a minorization condition for the $i$-times conditional filter in 
terms of the expected value of $\normalizerEst$ under the $(i+1)$-times conditional filter. 
The remainder of the section presents conditions under which the expected value of 
$\normalizerEst$ can be bounded. 
 One of the key assumptions is that adaptation controls the $\infty$-ESS.

Of particular interest are the cases $i=0$, which corresponds to the 
$\alpha$SMC filter, and $i=1$, which corresponds to
the conditional $\alpha$SMC kernel.
The former provides approximate samples from $\truedAll$.
The latter can be used to define 
a Markov chain with invariant distribution $\truedAll$, producing an adaptive resampling particle Gibbs sampler. 
We consider both these applications in, respectively, \cref{sec:divergence,sec:particle-mcmc}.

\subsection{A minorization condition for the \ciSMC kernel}

For the remainder of the article, we will work 
under the following assumption:
\begin{assumption} \label{asm:doubly-stochastic}
For all $N \ge 1$, all $\alpha \in \mathbb A_{N}$ are doubly stochastic.
\end{assumption}
\begin{remark} 
\cref{asm:doubly-stochastic} is the same as Assumption (B$^{++}$) in 
\citep{Whiteley:2013}, although there, the condition is stated as 
assuming each $\alpha$ admits the uniform distribution as
an invariant measure. 
\end{remark}

Let
\[
\kappa_{N} \defined \max_{\substack{n \ne m,\alpha \in \mathbb A_{N}}} \sum_{k=1}^{N}\alpha^{kn}\alpha^{km}
\qquad \text{and} \qquad
\kappa_{N}' \defined \max\{\kappa_{N}, 1/N\}.
\]
Our first main result provides control over how much the measure $\ciKernely$ differs 
from $\truedAll$. 
The theorem gives a stronger result when $i=0$ and 
gives a simpler result when $i=1$, by expressing the lower bound on 
$\ciKernel[1]{}(\yAll, S)$ in terms the more transparent quantity $\kappa_{N}'$ 
instead of the normalization terms $\ciNormConst{2}$. 
For example, if $\mathbb A_{N} = \{ \identityMatrix \}$ then $\kappa_{N} = 0$, while if 
$\mathbb A_{N} = \{ \uniformMatrix \}$, then $\kappa_{N} = 1/N$, so in either case 
$\kappa_{N}' = 1/N$. 

\bthm \label{thm:alpha-smc-rd-general}
If \cref{asm:doubly-stochastic} holds, 
then for all 
$t \ge 1$, 
$i \ge 0$, 
$N > i$, 
$S \subseteq \X^{t}$ measurable,
and $\yAll^{1},\dots,\yAll^{i},\xAll \in \X^{t}$, 
\[
\ciKernel{}(\yiAll{i}, S)
&\ge \int_{S}\frac{\normalizer}{\ciEE[i+1]{\yiAll{i},\xAll}[\normalizerEst \prod_{s=1}^{t}\ciNormConst{i+1}/\ciNormConst{i}]} \truedAll[t](\dee  \xAll).
\label{eq:ci-alpha-smc-bound}
\]
In particular, in the case of $i=0$, we have
\[
\frac{\dee \ciKernel[0]{}}{\dee \truedAll[t]}( \xAll)
&= \ciEE[1]{ \xAll}\left[\frac{\normalizer}{\normalizerEst}\right]
\ge \frac{\normalizer}{\ciEE[1]{ \xAll}[\normalizerEst]} 
\label{eq:alpha-smc-rd}
\]
and in the case of $i=1$, we have
\[
\ciKernel[1]{}(\yAll, S)
&\ge \int_{S}\frac{\normalizer(1-\kappa_{N}')^{t}}{\ciEE[2]{\yAll,\xAll}[\normalizerEst]} \truedAll[t](\dee  \xAll).
\label{eq:c-alpha-smc-bound}
\]
\ethm

\begin{remark}
By identical arguments, \cref{thm:alpha-smc-rd-general} also holds in the 
marginal and predictive cases.
In the predictive cases, however, $\normalizerEstPred$ and $\normalizerPred$ 
replace, respectively, $\normalizerEst$ and $\normalizer$.
In the predictive case, under \cref{asm:doubly-stochastic}, 
$\normalizerEstPred = \normalizerEst[t-1]$ and $\normalizerPred = \normalizer[t-1]$, 
so later results pertaining to ${\ciEEy}[\normalizerEst]$, such as 
\cref{prop:expected-sir-normalizer,prop:expected-smc-normalizer}, apply to $\normalizerPred$ as well.
The fact that $\normalizerPred = \normalizer[t-1]$ follows immediately from the definitions.
To show that $\normalizerEstPred = \normalizerEst[t-1]$, apply \cref{asm:doubly-stochastic}:
\begin{align*}
\normalizerEstPred[t] 
&= N^{-1} \sum_{n} \Wt{t}  
= N^{-1} \sum_{n} \sum_{k} \alpha_{t-1}^{nk}\Wt[k]{t-1} \potential[t-1]^{k} 
= N^{-1} \sum_{k} \Wt[k]{t-1} \potential[t-1]^{k}
= \normalizerEst[t-1].\qquad\quad \qedhere
\end{align*}
\end{remark}

We will prove \cref{thm:alpha-smc-rd-general} in two parts:
first for the case of $i = 0$, then for $i \ge 1$. 
For the $i=0$ case (corresponding to vanilla $\alpha$SMC), 
we begin by writing the joint density of the $\alpha$SMC process as
\[
\begin{split}
\lefteqn{\psi(\samps, \ancestors[t-1])} \\
&\defined \left(\prod_{n=1}^{N} \propd[1](\samp{1})\right)\left( \prod_{s=2}^{t} \prod_{n=1}^{N}r_{n}(\ancestor{s-1} | \wts[s-1], \samps[s-1])\propd[s](\samp[\ancestor{s-1}]{s-1}, \samp{s})\right).
\end{split}
\]
Under \cref{asm:doubly-stochastic}, the density of the \cSMC process with law 
$\cPry[\Samps, \Ancestors[t-1], \LineageFull[1]{t}]$ 
can be written in the following ``collapsed'' form, 
by implicitly identifying $\sampAll[{\lineage[1]{t}}]{t}$ with $\yAll[t]$:
\[
\lefteqn{\tpsi(\samps, \ancestors[t-1], \lineageFull[1]{t})} \nonumber \\
&= \frac{\psi(\samps, \ancestors[t-1]) \prod_{s=2}^{t}I_{s}\alpha_{s-1}^{\lineage[1]{s}\lineage[1]{s-1}}}{N\propd[1](\samp[{\lineage[1]{1}}]{1}) \prod_{s=2}^{t} r_{\lineage[1]{s}}(\lineage[1]{s-1} | \wts[s-1], \samps[s-1]) \propd[s](\samp[{\lineage[1]{s-1}}]{s-1}, \samp[{\lineage[1]{s}}]{s})} \\
&=  \frac{1}{N}\prod_{n \ne \lineage[1]{1}} \propd[1](\samp{1}) \prod_{s=2}^{t}\left(I_{s} \alpha_{s-1}^{\lineage[1]{s}\lineage[1]{s-1}} \prod_{n \ne \lineage[1]{s}} r_{n}(\ancestor{s-1} | \wts[s-1], \samps[s-1])\propd[s](\samp[\ancestor{s-1}]{s-1}, \samp{s})\right), \nonumber
\]
where $I_{s} \defined \ind(a_{s-1}^{\lineage[1]{s}} = \lineage[1]{s-1})$.

\bprf[Proof of \cref{thm:alpha-smc-rd-general}, $i = 0$ case]
Consider the density 
\[
\smcCondKernelDist(\samps, \ancestors[t-1], \lineageFull[1]{t}) 
&\defined \truedAll[t](\xAll[t]^{\lineage[1]{t}}) \tpsi(\samps, \ancestors[t-1], \lineageFull[1]{t}). \label{eq:artifical-target-measure}
\]
Then
\(
\lefteqn{\frac{\psi(\samps, \ancestors[t-1])\potential[t](\samp[{\lineage[1]{t}}]{t})\wt[{\lineage[1]{t}}]{t}}{\smcCondKernelDist(\samps, \ancestors[t-1], \lineageFull[1]{t})\sum_{n=1}^{N} \potential[t](\samp{t}) \wt{t}}}  \\
&= \frac{\propd[1](x_{1}^{\lineage[1]{1}}) \prod_{s=2}^{t} r_{\lineage[1]{s}}(\lineage[1]{s-1} | \wts[s-1], \samps[s-1]) \propd[s](\samp[{\lineage[1]{s-1}}]{s-1}, \samp[{\lineage[1]{s}}]{s}) \potential[t](\samp[{\lineage[1]{t}}]{t})\wt[{\lineage[1]{t}}]{t}}{\truedAll[t](\xAll[t]^{\lineage[1]{t}}) \prod_{s=2}^{t}I_{s}\alpha_{s-1}^{\lineage[1]{s}\lineage[1]{s-1}}N^{-1}\sum_{n=1}^{N}  \potential[t](\samp{t})\wt{t}} \\
&= \frac{\propd[1](x_{1}^{\lineage[1]{1}}) \prod_{s=2}^{t} \alpha_{s-1}^{\lineage[1]{s}\lineage[1]{s-1}}\wt[{\lineage[1]{s-1}}]{s-1}\potential[s-1](\samp[{\lineage[1]{s-1}}]{s-1}) \propd[s](\samp[{\lineage[1]{s-1}}]{s-1}, \samp[{\lineage[1]{s}}]{s})\potential[t](\samp[{\lineage[1]{t}}]{t})\wt[{\lineage[1]{t}}]{t}}{\truedAll[t](\xAll[t]^{\lineage[1]{t}}) \prod_{s=2}^{t} \wt[{\lineage[1]{s}}]{s} \prod_{s=2}^{t}I_{s}\alpha_{s-1}^{\lineage[1]{s}\lineage[1]{s-1}}N^{-1}\sum_{n=1}^{N} \potential[t](\samp{t})\wt{t}} \\
&= \frac{\prod_{s=1}^{t}\potential[s](\samp[{\lineage[1]{s}}]{s}) \propd[s](\samp[{\lineage[1]{s-1}}]{s-1}, \samp[{\lineage[1]{s}}]{s})}{\truedAll[t](\xAll[t]^{\lineage[1]{t}})N^{-1}\sum_{n=1}^{N} \potential[t](\samp{t})\wt{t}\prod_{s=2}^{t}I_{s}}  \\
&= \frac{\normalizer[t]}{\normalizerEst}\frac{1}{\prod_{s=2}^{t}I_{s}},
\)
Using the convention that $0/0 = 0$, it follows that
\(
\lefteqn{\ciKernel[0]{}(\dee \xAll)} \\
&= \sum_{\ancestors[t-1],\ancestorFinal} \int \Bigg\{\psi(\samps, \ancestors[t-1])\frac{\potential[t](\samp[\ancestorFinal]{t})\wt[\ancestorFinal]{t}}{\sum_{n=1}^{N} \potential[t](\samp{t})\wt{t}} \delta_{\xAll[t]^{a_t}}(\dee \xAll)\Bigg\}\dee\samps \\
\begin{split}
&= \sum_{\ancestors[t-1],\lineageFull[1]{t}} \int \Bigg\{\frac{\psi(\samps, \ancestors[t-1])\potential[t](\samp[{\lineage[1]{t}}]{t})\wt[{\lineage[1]{t}}]{t}\prod_{s=2}^{t}I_{s}}{\smcCondKernelDist(\samps, \ancestors[t-1], \lineages)\sum_{n=1}^{N}\potential[t](\samp{t}) \wt{t}} \\
&\phantom{\sum_{\ancestors[t-1],\lineages}~\int~\Bigg\{~} \times \smcCondKernelDist(\samps, \ancestors[t-1], \lineageFull[1]{t})\delta_{\xAll[t]^{f_t}}(\dee \xAll)\Bigg\}\dee\samps
\end{split} 
\\
&= \sum_{\ancestors[t-1],\lineageFull[1]{t}} \int \Bigg\{ \frac{\normalizer[t]}{\normalizerEst} \smcCondKernelDist(\samps, \ancestors[t-1], \lineageFull[1]{t})\delta_{\xAll[t]^{f_t}}(\dee \xAll)\Bigg\}\dee\samps \\
&= \Bigg\{ \sum_{\ancestors[t-1],\lineageFull[1]{t}} \int \frac{\normalizer[t]}{\normalizerEst} \tpsi(\samps, \ancestors[t-1], \lineageFull[1]{t})\delta_{\xAll}(\dee \xAll[t]^{f_t})\dee\bx_{1:t}^{-\lineage[1]{t}}\Bigg\} \truedAll[t](\dee \xAll).
\)
The result follows from \cref{lem:indicator}.
\eprf
We defer the proof of \cref{thm:alpha-smc-rd-general} in the $i \ge 1$ case to \cref{app:alpha-smc-rd-general-proof}.

\subsection{Bounding $\normalizerEst$ under the \ciSMC kernel}

In order to apply \cref{thm:alpha-smc-rd-general}, we must be able to 
control the quantity ${\ciEEy}[\normalizerEst]$.
As an initial step toward this goal, we consider the SIR case:

\setassumptiontag{SIR}
\begin{assumption} \label{asm:sir}
For all $s \in [t-1]$, $\alpha_{s} = \uniformMatrix$. 
\end{assumption}

For SIR, we can rewrite ${\ciEEy}[\normalizerEst]$ in an equivalent but more explicit 
form (\cref{prop:expected-sir-normalizer}).
Our goal will then be to provide general conditions under which ${\ciEEy}[\normalizerEst]$ can be rewritten 
in a similar manner (\cref{prop:expected-smc-normalizer}).

Recall that $G_{0,t} \defined \EE[\potentials[1]{t-1}(\all{\xi}{t-1})]$ and 
$G_{s,t}(x_{s}) \defined \EE[\potentials[s]{t-1}(\range{\xi}{s}{t-1}) \given \xi_{s} = x_{s}]$
for $s \in [t]$ and $x_{s} \in \X$.
For $t \ge 1$, $1 \le \ell \le s \le t + 1$, let
\[
\mcT_{t,\ell,s} \defined \theset{ \tuple[\tau_{1},\dots,\tau_{\ell}] : t - s + 1 < \tau_{1} < \dots < \tau_{\ell} = t + 1}
\]
and, for $\btau \in \mcT_{t,\ell,s}$, define
\[
C_{\ell}(\btau, \yAll[t]) \defined \prod_{i=1}^{\ell-1} G_{\tau_{i}, \tau_{i+1}}(y_{\tau_{i}}).
\]
We will sometimes write $C_{\ell}^{y}(\btau)$ or $C_{\ell}^{\btau}(\yAll[t])$
instead of $C_{\ell}(\btau, \yAll[t])$. 
The following is a straightforward generalization of \citep[Proposition 9]{Andrieu:2013}.

\bprop \label{prop:expected-sir-normalizer}
If \cref{asm:sir} holds, then 
for all $t \ge 1$, $i \ge 1$, $N \ge i$, $\yAll^{1},\dots,\yAll^{i} \in \X^{t}$, 
\[
{\ciEEy}[\normalizerEst]
= \frac{1}{N^{t}}\sum_{\ell=1}^{t+1}(N-i)^{t+1-\ell}
  \sum_{\btau \in \mcT_{t,\ell,t+1}} G_{0,\tau_{1}}  \prod_{m=1}^{\ell-1}\sum_{j=1}^{i} G_{\tau_{m}, \tau_{m+1}}(y_{\tau_{m}}^{j}).
\label{eq:expected-sir-normalizer}
\]
In particular, in the case of $i=1$, we have
\[
\ciEE[1]{\yAll}[\normalizerEst]
 = \frac{1}{N^{t}}\sum_{\ell=1}^{t+1}(N-1)^{t+1-\ell}
     \sum_{\btau \in \mcT_{t,\ell,t+1}} G_{0,\tau_{1}}C_{\ell}(\btau, \yAll[t]).
     \label{eq:sir-normalizer}
\]
\eprop

In order to obtain a version of \cref{prop:expected-sir-normalizer} for the general \ciSMC case,
we will require that the algorithm enforce a lower bound on a carefully chosen notion of effective
sample size called $\infty$-ESS. 
The $\infty$-ESS is a member of a family of effective sample size measures we call $p$-ESS, 
which also includes two commonly used definitions as special cases.

\bdefo
For parameter $p \in [1,\infty]$, let $p_{*} \defined \frac{p}{p - 1}$ be the
conjugate exponent of $p$ (so $1/p + 1/p_{*} = 1$). 
The \textbf{$p$-effective sample size} ($p$-ESS) of the weight vector 
$w^{1:N} \in \reals_{+}^{N}$ is
\[
\ess_{p}(w^{1:N}) \defined \left\{ 
  \begin{array}{l l}
    \left(\frac{\|w^{1:N}\|_{1}}{\|w^{1:N}\|_{p}}\right)^{p_{*}} & \quad p > 1 \\
    \frac{\|w^{1:N}\|_{1}}{\prod_{n=1}^{N} (w^{n})^{w^{n}/\|w^{1:N}\|_{1}}} & \quad p = 1.
  \end{array} \right.
\]
\edefo
The following proposition highlights some elementary 
properties of $p$-ESS and subsumes \cref{prop:basic-ess-properties} 
(see \cref{app:ess-properties-proof} for a proof).
\bprop \label{prop:ess-properties}
The $p$-ESS has the following properties:
\benum[leftmargin=*]
\item For all $p \in [1,\infty]$, $1 \le \ess_{p}(w^{1:N}) \le N$. 
The lower bound is achieved if and only if all but one of the weights is zero. 
The upper bound is achieved if and only if all the weights are equal.
\item For $1 < p < q \le \infty$, $\ess_{p}(w^{1:N}) \ge \ess_{q}(w^{1:N}) \ge N^{-(1 - q_{*}/p_{*})}\ess_{p}(w^{1:N})$,
with equality if and only if $K$ weights are equal and the rest are zero.  
\item The 1-ESS satisfies
\[
\ess_{1}(w^{1:N}) = \lim_{p \downarrow 1} \ess_{p}(w^{1:N}) = \ess_{ent}(w^{1:N}) \defined e^{H(w^{1:N})},
\]
where $H(w^{1:N}) \defined - \sum_{n} \frac{w^{n}}{\|w^{1:N}\|_{1}} \log \frac{w^{n}}{\|w^{1:N}\|_{1}}$
is the entropy. 
\eenum
\eprop

Parts (1) and (2) generalize their counterparts in 
\cref{prop:basic-ess-properties} to all $p \in [0,1]$, including the case $p=1$. 
Part (3) shows that the
$1$-ESS corresponds
to the entropic ESS, which is a common choice of 
ESS in applications \citep{Cornebise:2008}. 

In order to obtain a bound on ${\ciEEy}[\normalizerEst]$, we will require a lower
bound on the $\infty$-ESS of the weights, as formalized in \cref{asm:inf-ess-bound}.
Our development follows that of \citet{Whiteley:2013}, who used
the 2-ESS lower bound guarantee to bound the 
$L_{2}$ norm of the weights in terms of their $L_{1}$ norm.
Similarly, we will use the $\infty$-ESS lower bound
guarantee to bound the sup-norm of the weights in terms
of their $L_{1}$ norm.
Specifically, under \cref{asm:inf-ess-bound}, we have
\[
\zeta N 
&\le \ess_{\infty}(\Wts[s])
= \frac{\|\Wts[s]\|_{1}}{\|\Wts[s]\|_{\infty}} 
= \frac{\|\Wts[s]\|_{1}}{\sup_{n} \Wt{s}},
\]
and so, for all $n \in [N]$ and $s \in [t]$, 
$\Wt{s} \le \frac{\|\Wts[s]\|_{1}}{\zeta N}$.
We can use this upper bound on $\Wt{s}$ to prove a result that is very 
similar to \cref{prop:expected-sir-normalizer},
but permits an arbitrary adaptation scheme satisfying
\cref{asm:doubly-stochastic,asm:inf-ess-bound}.
(See \cref{app:expected-smc-normalizer-proof} for a proof):

\bthm \label{prop:expected-smc-normalizer}
If \cref{asm:doubly-stochastic,asm:inf-ess-bound} hold, 
then for all $t \ge 1$, $i \ge 1$, $N \ge i$, 
$\yAll^{1},\dots,\yAll^{i} \in \X^{t}$, 
\[
\lefteqn{{\ciEEy}[\normalizerEst]} \label{eq:expected-smc-normalizer-bound} \\
&\le \frac{1}{N(\zeta N)^{t-1}}\sum_{\ell=1}^{t+1} \sum_{\btau \in \mcT_{t,\ell,t+1}} (\zeta N)^{t+1-\ell}  \left(\frac{N - i}{\zeta N}\right)^{\ind(\tau_{1} > 1)}G_{0,\tau_{1}} \prod_{m=1}^{\ell-1}\sum_{j=1}^{i} G_{\tau_{m}, \tau_{m+1}}(y_{\tau_{m}}^{j}). \nonumber
\]
In particular, in the case of $i=1$, we have
\[
\ciEE[1]{\yAll}[\normalizerEst]
&\le \frac{1}{N(\zeta N)^{t-1}}\sum_{\ell=1}^{t+1} \sum_{\btau \in \mcT_{t,\ell,t+1}} (\zeta N)^{t+1-\ell}  \left(\frac{N - 1}{\zeta N}\right)^{\ind(\tau_{1} > 1)}G_{0,\tau_{1}}C_{\ell}^{y}(\btau).
\]
\ethm

The gap between \cref{prop:expected-sir-normalizer,prop:expected-smc-normalizer} is that
most of the factors of $N - i$ in the former are replaced by factors of $N$ in the latter.
Luckily we are interested in the $i=1,2$ cases, so we expect the differences between the 
two quantities to be fairly small. 
The following result, which is immediate upon expanding the left-hand sides of 
\cref{eq:expected-sir-normalizer,eq:expected-smc-normalizer-bound} and keeping only
$\Omega(1/N)$ terms, formalizes this intuition: 

\bcor \label{cor:expected-smc-normalizer-order}
If \cref{asm:sir} holds, then for all $i,t \ge 1$, $N > i$, 
$\yAll^{1},\dots,\yAll^{i} \in \X^{t}$,
\[
{\ciEEy}[\normalizerEst/\normalizer] = 1 + \frac{\normalizer^{-1}\sum_{s=1}^{t}\sum_{j=1}^{i}G_{0,s}G_{s,t+1}(\yAll^{i}) - ti}{ N}  + \Theta(N^{-2}).
\]
If \cref{asm:doubly-stochastic,asm:inf-ess-bound} hold, then for all $i,t \ge 1$, $N > i$, 
$\yAll^{1},\dots,\yAll^{i} \in \X^{t}$, 
\[
{\ciEEy}[\normalizerEst/\normalizer]
&\le 1 + \frac{\normalizer^{-1}\sum_{s=1}^{t}\sum_{j=1}^{i}G_{0,s}G_{s,t+1}(\yAll^{i}) - \zeta i}{\zeta N}  + \Theta(N^{-2}).
\]
\ecor

\subsection{Quantitative bounds}
\label{sec:quantitative-bounds}

Recall \cref{asm:bounded-potentials,asm:G0-ratio-bound}, either of which can be used
in conjunction with \cref{prop:expected-smc-normalizer} to obtain uniform, quantitative bounds
on ${\ciEEy}[\normalizerEst]$ by following the approach of \citet{Andrieu:2013}: 

\boundedpotentials* 

\Gzeroratiobound*

\cref{asm:G0-ratio-bound} is implied by a standard ``strong mixing'' condition
which is often employed in SMC analyses (e.g., \citep{DelMoral:2004,Whiteley:2013}).
See \citet{Andrieu:2013} for details.

\bprop \label{prop:normalizer-simple-bound}
If $\alpha_{s} = \uniformMatrix$ for $s \in [t-1]$ and \cref{asm:bounded-potentials} holds, then 
for all $t \ge 1$, $i \ge 1$, $N \ge i$, $\yAll^{1},\dots,\yAll^{i} \in \X^{t}$, 
\[
{\ciEEy}[\normalizerEst/\normalizer] 
\le 1 + \biggl[\normalizer[t]^{-1}\prod_{s=1}^{t}\potentialUB[s] - 1\biggr]\left[1 - \left(1 - \frac{i}{N}\right)^{t}\right].
\label{eq:sir-normalizer-simple-bound}
\]
If \cref{asm:doubly-stochastic,asm:inf-ess-bound,asm:bounded-potentials} hold, then 
for all $t \ge 1$, $i \ge 1$, $N \ge i$, $\yAll^{1},\dots,\yAll^{i} \in \X^{t}$, 
\[
{\ciEEy}[\normalizerEst/\normalizer] 
\le 1 + \normalizer[t]^{-1}\prod_{s=1}^{t}\potentialUB[s]\left[\left(1 + \frac{i}{\zeta N}\right)^{t} - 1\right].
\label{eq:smc-normalizer-simple-bound}
\]
\eprop

\bprf
The proof \cref{eq:sir-normalizer-simple-bound} is a straightforward generalization of that for 
\citep[Proposition 12]{Andrieu:2013} with some additional bookkeeping for  
$i$ (instead of 2) fixed trajectories.
As for \cref{eq:smc-normalizer-simple-bound}, we have 
\[
\begin{split}
{\ciEEy}[\normalizerEst[t]] 
&\le \sum_{\ell=1}^{t+1} \sum_{\btau \in \mcT_{t,\ell,t+1}} (\zeta N)^{-\ell+1}G_{0,\tau_{1}}\prod_{m=1}^{\ell-1}\sum_{j=1}^{i} G_{\tau_{m}, \tau_{m+1}}(y_{\tau_{m}}^{j}) \\
&\le \normalizer + \prod_{s=1}^{t}\potentialUB[s]\sum_{\ell=2}^{t+1} {t \choose \ell-1} i^{\ell-1} (\zeta N)^{-\ell+1}  \\
&= \normalizer + \prod_{s=1}^{t}\potentialUB[s]\sum_{\ell=1}^{t} {t \choose \ell} (\zeta N/i)^{-\ell}  \\
&= \normalizer + \prod_{s=1}^{t}\potentialUB[s]\left[\left(1 + \frac{i}{\zeta N}\right)^{t} - 1\right] .
\end{split}
\]
\eprf

\bprop \label{prop:normalizer-G0-bound}
If $\alpha_{s} = \uniformMatrix$ for $s \in [t-1]$ and \cref{asm:G0-ratio-bound} holds, then 
for all $t \ge 1$, $i \ge 1$, $N \ge i$, $\yAll^{1},\dots,\yAll^{i} \in \X^{t}$, 
\[
{\ciEEy}[\normalizerEst/\normalizer]
\le \left(1 + \frac{i(\beta-1)}{N}\right)^{t}.
\label{eq:sir-normalizer-G0-bound}
\]
If \cref{asm:doubly-stochastic,asm:inf-ess-bound,asm:G0-ratio-bound} hold, then 
for all $t \ge 1$, $i \ge 1$, $N \ge i$, $\yAll^{1},\dots,\yAll^{i} \in \X^{t}$, 
\[
{\ciEEy}[\normalizerEst/\normalizer] 
\le  \left(1 + \frac{i\beta}{\zeta N}\right)^{t}.
\label{eq:smc-normalizer-G0-bound}
\]
\eprop

\bprf
The proof of \cref{eq:sir-normalizer-G0-bound} is a simple generalization of that for 
\citep[Proposition 14]{Andrieu:2013}.
As for \cref{eq:smc-normalizer-G0-bound}, observe that for $s \in [t+1]$, 
$G_{0,t+1} = G_{0,t+1}\frac{G_{0,t+1}}{G_{0,s}} = G_{0,s}\truedPredMarginal[s](G_{s,t+1})$,
so we can write for $\ell \in [t]$, $\btau \in \mcT_{t,\ell,t+1}$,
\[
\normalizer = G_{0,t+1} 
= G_{0,\tau_{k}}\prod_{i=1}^{\ell-1} \truedPredMarginal[\tau_{i}](G_{\tau_{i},\tau_{i+1}}).
\]
Combined with \cref{asm:G0-ratio-bound} and writing
$\bar G_{s,t} \defined \sup_{x \in \X} G_{s,t}(x)$,
\[
\begin{split}
\lefteqn{\sum_{\ell=1}^{t+1} (\zeta N)^{-\ell+1}\sum_{\tau \in \mcT_{t,\ell,t+1}}G_{0,\tau_{1}}\prod_{i=1}^{\ell-1}\sum_{j=1}^{i} G_{\tau_{m}, \tau_{m+1}}(y_{\tau_{m}}^{j})} \\
&\le \normalizer + \normalizer \sum_{\ell = 2}^{t+1} (\zeta N)^{-\ell+1} \sum_{\tau \in \mcT_{t,\ell,t+1}}\frac{G_{0,\tau_{1}}}{G_{0,\tau_{1}}} \prod_{i=1}^{\ell-1} \sum_{j=1}^{i}\frac{\bar G_{\tau_{i},\tau_{i+1}}}{\truedPredMarginal[\tau_{i}](G_{\tau_{i},\tau_{i+1}})} \\
&= \normalizer \sum_{\ell=1}^{t+1}{t \choose \ell - 1}(\zeta N)^{-\ell+1} (i\beta)^{\ell-1} \\
&= \normalizer \left(1 + \frac{i\beta}{\zeta N}\right)^{t}.
\end{split}
\]
\eprf

To compare \cref{eq:sir-normalizer-simple-bound,eq:smc-normalizer-simple-bound}, consider the $\Theta(1/N)$
terms, which are, respectively, 
\(
\frac{ti[\normalizer[t]^{-1}\prod_{s=1}^{t}\potentialUB[s] - 1]}{N}
\qquad \text{and} \qquad
\frac{ti\normalizer[t]^{-1}\prod_{s=1}^{t}\potentialUB[s]}{\zeta N}.
\)
Thus, up to a $-ti$ term and a factor of $1/\zeta$, the two bounds are of the same
leading order in $1/N$. 
The $-it$ is likely an artifact of the analysis while the $1/\zeta$ term accounts for
there being only $\zeta N$ ``effective particles.''
The differences between \cref{eq:sir-normalizer-G0-bound,eq:smc-normalizer-G0-bound}
are identical.

\section{Bounding the Divergence of SMC Samplers}
\label{sec:divergence}

Recall that $\ciKernel[0]{}(\dee \xAll)$ is the distribution of $\SampAllPathFinal \dist \truedEstT[0]{t}$,
a single sample from the $\alpha$SMC estimator of $\truedAll$. 
Equivalently, 
$\ciKernel[0]{}(\dee \xAll) = \ciEE[0]{}[\truedEstT[0]{t}](\dee \xAll) \defined \ciEE[0]{}[\truedEstT[0]{t}(\dee \xAll)]$
is the expected value of the random measure $\truedEstT[0]{t}$.
As a first application of our results from \cref{sec:main}, we consider
bounding the distance between the measures $\truedAll$ and $\ciKernel[0]{}$.
That is, for some divergence $d(\mu || \nu)$ between measures, can we bound 
$d(\truedAll || \ciKernel[0]{})$?
To the best of our knowledge, there has been minimal investigation of this question, with
\citep[Chapter 8]{DelMoral:2004} a notable exception.
For example, under \cref{asm:sir}, the bound
\[
\kl{\ciKernel[0]{}}{\truedMarginal} \le \frac{c}{N}, \label{eq:reverse-kl-bound}
\]
can be extracted as a special case of a more general propagation-of-chaos 
result \citep[Theorem 8.3.2]{DelMoral:2004}. 

\newcommand{\MID}{\mathcal F_1}
Let $\MID$ be the set of functions
$f : \reals_{+} \to \reals$ that are monotonically increasing or decreasing and satisfy $f(1) = 0$.
We consider the class of \emph{monotonic divergences} of the form
\[
d_{i,f}(\mu_{1} || \mu_{2}) 
\defined \mu_{i} \Bigl (f \circ \frac{\dee \mu_{1}}{\dee \mu_{2}} \Bigr), 
\quad i \in \{1,2\}, f \in \MID. 
\label{eq:monotonic-divergence}
\]
\cref{tbl:monotonic-divergences} lists some common divergences that can be written this
way.

The following result characterizes the divergence between $\truedAll$ and $\ciKernel[0]{}$, only assuming that $f$ is concave.

\bprop \label{prop:smc-sampler-divergences}
Let $\mcR_{t}(\yAll) \defined {\cEEy}[\normalizerEst/\normalizer]$ and 
$\mcS_{t} \defined \normalizer^{-1}\sum_{s=1}^{t}G_{0,s}\truedAll(G_{s,t+1}) - \zeta$. 
If \cref{asm:doubly-stochastic} holds, then for all concave $f \in \MID$,
\(
d_{1,f}(\truedAll || \ciKernel[0]{}) \le f(\truedAll(\mcR_{t})).
\)
In particular, if \cref{asm:sir} holds, then 
\(
\kl{\truedAll}{\ciKernel[0]{}} &\le \log\left(1 + \frac{\mcS_{t} - t}{N} + \Theta(N^{-2})\right) \\
\chisq{\truedAll}{\ciKernel[0]{}} &\le \frac{\mcS_{t} - t}{N} + \Theta(N^{-2})
\)
while if \cref{asm:doubly-stochastic,asm:inf-ess-bound} hold, then
\(
\kl{\truedAll}{\ciKernel[0]{}} &\le \log\left(1 + \frac{\mcS_{t} - \zeta}{\zeta N} + \Theta(N^{-2})\right) \\
\chisq{\truedAll}{\ciKernel[0]{}} &\le \frac{\mcS_{t} - \zeta}{\zeta N} + \Theta(N^{-2}). 
\)
\eprop
\bprf
The general statement follows by applying Jensen's inequality 
to \cref{eq:monotonic-divergence}, then using \cref{thm:alpha-smc-rd-general}.
The special cases correspond to using the KL divergence (version 1) and $\chi^{2}$ distance
rows of \cref{tbl:monotonic-divergences} and applying  \cref{cor:expected-smc-normalizer-order}.
\eprf

Similar results for (sequential) importance sampling are included in \cref{app:is}.

We can also consider the divergence between $\truedAll$ and $\ciKernel[0]{}$
when ${\cEEy}[\normalizerEst]$ is uniformly bounded:

\bprop \label{prop:smc-sampler-divergences}
If \cref{asm:doubly-stochastic} holds and ${\cEEy}[\normalizerEst/\normalizer] \le B_{t,N}$
for all $\yAll[t] \in \X^{t}$, then for all increasing $f \in \MID$
\(
d_{i,f}(\truedAll || \ciKernel[0]{}) \le f(B_{t,N})
\)
and for all decreasing $f \in \MID$,
\(
d_{i,f}(\ciKernel[0]{} || \truedAll) \le f(B_{t,N}^{-1}). 
\)
In particular,
\(
\kl{\truedAll}{\ciKernel[0]{}} &\le \log B_{t,N} \\
\chisq{\truedAll}{\ciKernel[0]{}} &\le B_{t,N} - 1 \\
\totalvar{\truedAll}{\ciKernel[0]{}} &\le \frac{B_{t,N} - 1}{B_{t,N}} \le B_{t,N} - 1.
\)
\eprop
\bprf
The general statements follow immediately from \cref{eq:monotonic-divergence,thm:alpha-smc-rd-general}. 
The special cases correspond to using the KL divergence (version 1), $\chi^{2}$ distance,
and total variation distance (version 2) rows of \cref{tbl:monotonic-divergences}.
The second total variation inequality holds since $B_{t,N} \ge 1$. 
\eprf

\begin{table}[t]
\caption{Divergences of the form \cref{eq:monotonic-divergence}.
The operator $(a)^{+}$ gives the positive part of $a \in \reals$.}
\begin{center}
\begin{tabular}{l|l|c|c}
Name & Symbol & $i$ & $f$ \\
\hline
KL divergence (version 1) & $d_{KL}$ & 1 & $a \mapsto \log a$ \\
KL divergence (version 2) & $d_{KL}$ & 2 & $a \mapsto - \log a$ \\
$\chi^{2}$ distance & $d_{\chi^{2}}$ & 1 & $a \mapsto a - 1$ \\
total variation distance (version 1) & $d_{TV}$ & 2 & $a \mapsto (a - 1)^{+}$ \\
total variation distance (version 2) & $d_{TV}$ & 2 & $a \mapsto (1 - a)^{+}$ \\
\end{tabular}
\end{center}
\label{tbl:monotonic-divergences}
\end{table}

The bounds in \cref{prop:smc-sampler-divergences} for KL divergence, $\chi^{2}$ distance,  and total variation distance are asymptotically
equivalent if $B_{t,N} \to 1$ as $N \to \infty$. 
Combining \cref{prop:smc-sampler-divergences} with, for example, 
\cref{prop:normalizer-G0-bound}, yields quantitative bounds for
SIR and $\alpha$SMC:

\bcor
If \cref{asm:sir,asm:G0-ratio-bound} hold, then 
\[
\kl{\truedAll}{\ciKernel[0]{}} &\le \frac{t(\beta-1)}{N} \\
\chisq{\truedAll}{\ciKernel[0]{}} &\le \frac{t(\beta-1)}{N} + O(N^{-2}) \\
\totalvar{\truedAll}{\ciKernel[0]{}} &\le \frac{t(\beta-1)}{N + t(\beta-1)} + O(N^{-2}).
\]
If \cref{asm:doubly-stochastic,asm:inf-ess-bound,asm:G0-ratio-bound} hold, then
\[
\kl{\truedAll}{\ciKernel[0]{}} &\le \frac{t\beta}{\zeta N} \\
\chisq{\truedAll}{\ciKernel[0]{}} &\le \frac{t\beta}{\zeta N} + O(N^{-2}) \\
\totalvar{\truedAll}{\ciKernel[0]{}} &\le \frac{t\beta}{\zeta N + t\beta} + O(N^{-2}).
\]
\ecor

\section{The $\alpha$-Particle Gibbs Sampler}
\label{sec:particle-mcmc}

As a second application of our results from \cref{sec:main}, we 
consider the mixing properties of particle Gibbs with adaptive resampling.
Recall from \cref{sec:intro-applications} that we introduce a global parameter
$\theta \in \Theta$ with prior distribution $\varpi(\dee\theta)$.
Replace $\propd[s]$ by $\propdParam{s}$ and 
$\potential[s]$ by $\potentialParam{s}$, then parameterize
the other quantities defined previously in terms of $\propd[s]$ 
and $\potential[s]$ by $\theta$.
The target distribution on the product space 
$(\Theta \times \Xtsh, \mcB(\Theta \times \Xtsh))$, $\Xtsh \defined \X^{t}$, is 
\[
\trued(\dee \theta \times \dee y)
&\defined \truedU(\dee \theta \times \dee y) / \normalizer[],
\]
where
\[
\truedU(\dee \theta \times \dee y) 
\defined \prod_{s=1}^{t} \potentialParam{s}(y_{s})\propdParam{s}(y_{s-1}, \dee y_{s}) \varpi(\dee \theta)
\quad \text{and} \quad
\normalizer[] \defined \truedU(1). 
\]

Particle Gibbs samplers have kernels of the form
$\trued[y](\dee \theta)\Pi_{\theta}(y, \dee z)$,
where $\Pi_{\theta}(y, \dee z)$ is an SMC-based kernel
with invariant distribution $\trued[\theta]$. 
The standard PG sampler employs the iterated conditional 
SMC (i-cSMC) kernel~\citep{Andrieu:2013}:
that is, $\Pi_{\theta} = \cKernel[\theta]$ and require
\cref{asm:sir} to hold.

We now introduce the novel $\alpha$-particle Gibbs ($\alpha$PG) sampler,
which employs the iterated conditional
$\alpha$SMC (i-\cSMC) kernel $\cKernel[\theta]$, so
$\Pi_{\theta} = \cKernel[\theta]$.
In \cref{sec:i-csmc-kernel-reversible} we prove %
that the i-\cSMC kernel is reversible with respect 
to $\trued[\theta]$ and hence has invariant distribution $\trued[\theta]$. 

The first step to proving mixing results for the i-\cSMC kernel
and the $\alpha$PG sampler is to use \cref{thm:alpha-smc-rd-general} 
to obtain a sufficient condition for the i-\cSMC transition kernel to 
satisfy a minorization condition.

\bprop \label{prop:i-c-smc-minorization}
If \cref{asm:doubly-stochastic} holds and ${\ccEE}[\normalizerEst/\normalizer] \le B_{t,N}$ 
for all $\theta \in \Theta$ and $y^{1}, y^{2} \in \Xtsh$, then
\[
\cKernel[\theta](y, \dee x) \ge \veps_{t,N} \trued[\theta](\dee x), \label{eq:minorization-condition}
\]
where $\veps_{t,N} \defined \frac{(1 - \kappa_{N}')^{t}}{B_{t,N}}$. 
\eprop

The constant $\eps_{t,N}$, which determines mixing speed, can be found explicitly 
using the quantitative bounds from \cref{sec:quantitative-bounds}.
For example, using \cref{asm:G0-ratio-bound} we obtain the following: 

\bcor \label{cor:i-c-smc-strong-bound}
If \cref{asm:doubly-stochastic,asm:inf-ess-bound,asm:G0-ratio-bound}
hold, then for all $y \in \Xtsh$, 
\[ 
\cKernel[\theta](y, \dee x) \ge \veps_{t,N}\trued[\theta](\dee x), 
\]
where 
\[
\veps_{t,N} \defined \biggl(\frac{1-\kappa_{N}'}{1 + \frac{2\beta}{\zeta N}} \biggr)^{t}.
\]
Furthermore, if 
$\zeta N \ge \frac{2\beta t}{C(1 - \kappa_{N}') - \kappa_{N}'t}$ 
for some constant $C > 0$, then 
\[
\eps_{t,N} \ge e^{-C}.
\]
In particular, assuming $\kappa_{N}' \le B/N$ for some constant
$B \ge 1$, if $N \ge t/C + B$, then
\[
\eps_{t,N} \ge \exp \Bigl(-\frac{C}{\zeta}(2\beta + \zeta B) \Bigr).
\]
\ecor
\bprf
The first part follows from \cref{prop:i-c-smc-minorization,prop:normalizer-G0-bound}.
For the second part, we then have 
\[
\eps_{t,N} 
&\ge \biggl(\frac{1 + \frac{2\beta}{\zeta N}}{1-\kappa_{N}'}\biggr)^{-t} 
= \biggl(1 + \frac{1}{1-\kappa_{N}'}\left(\frac{2\beta}{\zeta N} + \kappa_{N}'\right)\biggr)^{-t} \\
&\ge \biggl(1 + \frac{C}{t}\biggr)^{-t}
\ge e^{-C}. 
\]
The final part follows after noting that if
$\kappa_{N}' \le B/N$, then 
\[
\frac{1}{1-\kappa_{N}'}\left(\frac{2\beta}{\zeta N} + \kappa_{N}'\right)
\le \frac{1}{1-B/N}\left(\frac{2\beta}{\zeta N} + B/N\right)
= \frac{1}{N - B}\left(\frac{2\beta}{\zeta} + B\right).
\]
\eprf
\begin{remark}
In the case of the i-cSMC kernel, \cref{cor:i-c-smc-strong-bound}
is almost as tight as \citep[Corollary 14]{Andrieu:2013}: 
the former result replaces $\beta - 1$ with $\beta$. 
\end{remark}

The minorization condition \cref{eq:minorization-condition} implies uniform
ergodicity and a number of other types of convergence guarantees
for the i-cSMC process.
The following generalizes \citep[Theorem 1]{Andrieu:2013}, which applies
only to the i-cSMC kernel and the PG sampler.
\bthm
Assume that \cref{asm:doubly-stochastic,asm:inf-ess-bound} hold.
\benum[label=\Roman*.]
\item
Let $N \ge 2$, and consider the i-\cSMC process
with kernel $P = \ciKernel[1]{\theta}$. 
\benum[leftmargin=*,label=\arabic*.]
\item $P$ is reversible with respect to $\trued$ and defines
a positive operator,
\item If the potentials are bounded then there exists 
$\eps_{t,N} = 1 - C_{t}/N$ such that
\benum
\item for all $y \in \Xtsh$, 
$P(y, \dee z) \ge \eps_{t,N}\trued[\theta](\dee z)$,
\item for every measure $\nu \ll \trued[\theta]$ and $k \ge 1$,
\[
\chisq{\nu P^{k}}{\trued[\theta]}
\le \chisq{\nu}{\trued[\theta]}(1 - \eps_{t,N})^{k},
\]
\item for every $y \in \Xtsh$ and $k \ge 1$, 
\[
\totalvar{\delta_{y}P^{k}}{\trued[\theta]} \le (1 - \eps_{t,N})^{k},
\]
\eenum
\item If \cref{asm:G0-ratio-bound} also holds
and there is a constant $B > 0$ such that $\kappa_{N}' \le B/N$,
then for every $C > 0$, there exists $\veps_{B,C,\zeta} > 0$ such
that for $N \ge t/C + B$ and all $t > 1$,
\[
\eps_{t,N} \ge \veps_{B,C,\zeta} > 0.
\]
\eenum
\item
If there exists $\beta \ge 1$ such that, 
for all $t, s \in \naturals$,
\[
\trued\text{-}\esssup_{\theta,x} \frac{G_{0,t}^{\theta}G_{t,t+s}^{\theta}(x)}{G_{0,t+s}^{\theta}} \le \beta,
\]
or if 
\[
\trued\text{-}\esssup_{\theta} \frac{\prod_{s=1}^{t}\potentialUB[s]^{\theta}}{\truedU[\theta](1)} < \infty,
\]
then the $\alpha$PG chain is geometrically ergodic whenever the Gibbs sampler is geometrically 
ergodic.
\eenum
\ethm
\bprf
Part I.1 follows from \cref{lem:i-csmc-kernel-reversible}.
Parts I.2-3 follow from 
\cref{prop:i-c-smc-minorization}, \cref{cor:i-c-smc-strong-bound}, and
\citep[Proposition 31]{Andrieu:2013}.
Part II follows from \citep[Section 7]{Andrieu:2013}.
\eprf

\begin{remark}
Part I.3 means that if \cref{asm:G0-ratio-bound} holds, then scaling $N$ linearly
with $t$ ensures a uniform convergence rate, as measured by
$\chi^{2}$-divergence or total variation distance.
\end{remark}

\appendix

\section{Additional Proofs}

\subsection{Technical Lemma}

\blem \label{lem:indicator}
Let $X$ and $Y$ be random elements in Borel spaces $(S,\mcS)$ and $(T,\mcT)$, respectively, let $\psi : S \times T \to \nnreals$ be a measurable, and let $\mu$ be the distribution of $X$.
If 
\[
\nu = \EE [\psi(X,Y) \delta_{X} ],
\]
then $\nu \ll \mu$ and 
\[
\frac{\dee \nu}{\dee \mu}(X) = \EE [\psi(X,Y) \given X]\ \as
\]
\elem
\bprf 
Because $S$ is Borel, there exists an $f$ satisfying $f(X) = \EE [\psi(X,Y) \given X]$ a.s.  
It follows from the chain rule of conditional expectation and then some elementary manipulations that, for all $A \in \mcS$, 
\(
\nu(A) 
= \EE [ f(X) \delta_X(A)] 
= \EE [ f(X) \ind_A(X)]
= \int_{A} f(x) \mu(\dee x),
\)
and so $f$ is a version of the Radon--Nikodym derivative $\dee \nu/\dee \mu$. 
\eprf

\subsection{Proof of \cref{thm:alpha-smc-rd-general}, $i \ge 1$ case}
\label{app:alpha-smc-rd-general-proof}

First observe that we can write the \ciSMC kernel as
\[
\ciKernel{}(\yiAll{i}, \dee \xAll) 
&= \ciEEy\left[\delta_{\SampAllPathFinal}(\dee \xAll) \right] \\
&= \ciEEy\left[\textstyle\sum_{\all{k}{t} \in [N]^{t}} \mcI_{\all{k}{t}}(\Samps, \Ancestors, \dee \xAll)\right],
\]
where 
\[
\mcI_{\all{k}{t}}(\samps, \ancestors, \dee \xAll)
&\defined \delta_{\sampAllPath[k_{t}]{t}}(\dee \xAll) \ind(k_{t} = \ancestorFinal) \prod_{s=1}^{t-1} \ind(k_{s} = \ancestor[k_{s+1}]{s}).
\]
Next note that 
\[
\lefteqn{\sum_{k_{1}=1}^{N} \ind(\sampAllPath[k_{t}]{t} \in S) {\ciPry}[\SampsT[1] \in \dee\sampsT[1], \LineagesT[1] = \lineagesT[1]]} \\
&= \sum_{k_{1}=1}^{N} \ind(\sampAllPath[k_{t}]{t} \in S) \ciNormConst[1]{i} \mcD(\lineagesT[1])\prod_{j=1}^{i}\frac{1}{N}\delta_{y_{1}^{j}}(\dee \samp[{\lineage[j]{1}}]{1})
  \prod_{n \notin \lineagesT[1]}^{N} \propd[1](\dee \samp{1})  \\
\begin{split} 
&\ge \frac{N\,\ciNormConst[1]{i}}{\ciNormConst[1]{i+1}} \sum_{k_{1}=1}^{N} \int_{\X}\ind(\sampAllPath[k_{t}]{t} \in S) 
  \ciNormConst[1]{i+1}\mcD(\lineagesT[1], k_{1})\prod_{j=1}^{i}\frac{1}{N}\delta_{y_{1}^{j}}(\dee \samp[{\lineage[j]{1}}]{1})
  \frac{1}{N}\delta_{x_{1}}(\dee \samp[k_{1}]{1}) \\
&\phantom{\ge~\frac{N\,\ciNormConst[1]{i}}{\ciNormConst[1]{i+1}} \sum_{k_{1}=1}^{N} \int_{\X}}
  	 \times \prod_{n \notin \lineagesT[1], k_{1}}^{N} \propd[1](\dee \samp{1})\propd[1](\dee x_{1})  \\
\end{split} \\
&\ge \frac{N\,\ciNormConst[1]{i}}{\ciNormConst[1]{i+1}} \sum_{k_{1}=1}^{N} \int_{\X}\ind(\sampAllPath[k_{t}]{t} \in S) 
  \ciPr[i+1]{\yiAll{i},\xAll}[\SampsT[1] \in \dee\sampsT[1], \LineagesT[1] = \lineagesT[1], \Lineage[i+1]{1} = k_{1}] \propd[1](\dee x_{1}) \label{eq:sum-k-1}
\]
For the remainder of the proof, to keep notation compact
when writing laws, instead of writing, e.g.,
$\SampsT[s] \in \sampsT[s]$ or $\LineagesT[s] = \lineagesT[s]$, whenever
a random variable is instantiated to be (the differential) of
the lowercase version of itself, we will write only the 
random variable: for example, $\SampsT[s]$ or $\LineagesT[s]$. 
Now, for $s = 2,\dots,t$,
\[
\lefteqn{\sum_{k_{s}=1}^{N} \ind(\sampAllPath[k_{t}]{t} \in S) \ind(k_{s-1} = a_{s-1}^{k_{s}}){\ciPry}[\SampsT[s], \AncestorsT[s-1], \LineagesT[s] \given \Samps[s-1], \Ancestors[s-2], \LineagesT[s-1]]}  \\
\begin{split}
  &\defined \sum_{k_{s}=1}^{N} \ind(\sampAllPath[k_{t}]{t} \in S) \ind(k_{s-1} = a_{s-1}^{k_{s}})\ciNormConst{i} \mcD(\lineagesT[s]) \prod_{j=1}^{i}\alpha_{s-1}^{\lineage[j]{s}\lineage[j]{s-1}} \delta_{y_{s}^{j}}(\dee \samp[{\lineage[j]{s}}]{s})\ind(\ancestor[{\lineage[j]{s}}]{s-1} = \lineage[j]{s-1}) \\
  &\phantom{\defined~\sum_{k_{s}=1}^{N}} 
  	\times \prod_{n \notin \lineagesT[s]} r_{n}(\ancestor{s-1} | \wts[s-1], \samps[s-1])\propd[s](\samp[\ancestor{s-1}]{s-1}, \samp{s})
\end{split} \\
\begin{split}
&\ge  \sum_{k_{s}=1}^{N} \frac{\ciNormConst{i}}{\ciNormConst{i+1}\alpha_{s-1}^{k_{s}k_{s-1}}} \int_{\X} \ind(\sampAllPath[k_{t}]{t} \in S) \ciNormConst{i+1} \mcD(\lineagesT[s], k_{s}) \\
 &\phantom{=~\sum_{k_{s}=1}^{N}~\frac{\ciNormConst{i}}{\ciNormConst{i+1}\alpha_{s-1}^{k_{s}k_{s-1}}}\int} 
 	\prod_{j=1}^{i}\alpha_{s-1}^{\lineage[j]{s}\lineage[j]{s-1}} \delta_{y_{s}^{j}}(\dee \samp[{\lineage[j]{s}}]{s})\ind(\ancestor[{\lineage[j]{s}}]{s-1} = \lineage[j]{s-1})  \\
 &\phantom{=~\sum_{k_{s}=1}^{N}~\frac{\ciNormConst{i}}{\ciNormConst{i+1}\alpha_{s-1}^{k_{s}k_{s-1}}}\int}
 	\times \alpha_{s-1}^{k_{s}k_{s-1}}\delta_{x_{s}}(\dee \samp[k_{s}]{s}) \ind(k_{s-1} = a_{s-1}^{k_{s}}) \\
 &\phantom{=~\sum_{k_{s}=1}^{N}~\frac{\ciNormConst{i}}{\ciNormConst{i+1}\alpha_{s-1}^{k_{s}k_{s-1}}}\int}
 	 \times \prod_{n \notin \lineagesT[s], k_{s}} r_{n}(\ancestor{s-1} | \wts[s-1], \samps[s-1])\propd[s](\samp[\ancestor{s-1}]{s-1}, \samp{s}) \\
 &\phantom{=~\sum_{k_{s}=1}^{N}~\frac{\ciNormConst{i}}{\ciNormConst{i+1}\alpha_{s-1}^{k_{s}k_{s-1}}}\int}
 	 \times r_{k_{s}}(k_{s-1} | \wts[s-1], \samps[s-1]) \propd[s](\samp[k_{s-1}]{s-1}, \dee x_{s}) 
\end{split} \\
&= \sum_{k_{s}=1}^{N} \frac{\ciNormConst{i}}{\ciNormConst{i+1}\alpha_{s-1}^{k_{s}k_{s-1}}} \int_{\X} \ind(\sampAllPath[k_{t}]{t} \in S) r_{k_{s}}(k_{s-1} | \wts[s-1], \samps[s-1]) \propd[s](\samp[k_{s-1}]{s-1}, \dee x_{s}) \\
&\phantom{=~\sum_{k_{s}=1}^{N}}\ciPr[i+1]{\yiAll{i},\xAll}[\SampsT[s], \AncestorsT[s-1], \LineagesT[s], \Lineage[i+1]{s} = k_{s} \given \Samps[s-1], \Ancestors[s-2], \LineagesT[s-1], \Lineage[i+1]{s-1} = k_{s-1}] \nonumber
\label{eq:sum-k-s}
\]

Using \cref{eq:sum-k-1,eq:sum-k-s}, we have (note that the
terms such as those involving $a_{0}$ should be ignored)
\[
\lefteqn{\sum_{\all{k}{t} \in [N]^{t}} \ind(\sampAllPath[k_{t}]{t} \in S) \ind(k_{t} = \ancestorFinal) \prod_{s=2}^{t} \ind(k_{s-1} = a_{s-1}^{k_{s}}) {\ciPry}[\Samps, \Ancestors, \Lineages]} \nonumber \\
\begin{split}
&= \sum_{\all{k}{t} \in [N]^{t}} \ind(\sampAllPath[k_{t}]{t} \in S)\ind(k_{t} = \ancestorFinal) {\ciPry}[\AncestorFinal \given \Samps, \Ancestors[t-1]] {\ciPry}[\SampsT[1], \LineagesT[1]] \\
&\phantom{=~\sum~=} \times \prod_{s=2}^{t}\ind(k_{s-1} = a_{s-1}^{k_{s}}) {\ciPry}[\SampsT[s], \AncestorsT[s-1], \LineagesT[s] \given \Samps[s-1], \Ancestors[s-2], \LineagesT[s-1]]
\end{split} \\
&\ge \sum_{\all{k}{t} \in [N]^{t}} \int_{\X^{t}} \frac{N\,\ind(\sampAllPath[k_{t}]{t} \in S)\ind(k_{t} = \ancestorFinal)}{(\prod_{s=1}^{t} \ciNormConst{i+1}/\ciNormConst{i})(\prod_{s=1}^{t-1}\alpha_{s-1}^{k_{s}k_{s-1}})}  \frac{\wt[\ancestorFinal]{t}\potential(\samp[\ancestorFinal]{t})}{\sum_{n=1}^{N} \wt{t}\potential(\samp{t})} \nonumber \\
&\phantom{=N~\sum~} \times \prod_{s=2}^{t} r_{k_{s}}(k_{s-1} | \wts[s-1], \samps[s-1]) \propd[s](\samp[k_{s-1}]{s-1}, \dee x_{s}) \\
&\phantom{=N~\sum~} \times \prod_{s=2}^{t} \ciPr[i+1]{\yiAll{i},\xAll}[\SampsT[s], \AncestorsT[s-1], \LineagesT[s], \Lineage[i+1]{s} = k_{s} \given \Samps[s-1], \Ancestors[s-2], \LineagesT[s-1], \Lineage[i+1]{s-1} = k_{s-1}] \nonumber \\
&\phantom{=N~\sum~} \times  \ciPr[i+1]{\yiAll{i},\xAll}[\SampsT[1] \in \dee\sampsT[1], \LineagesT[1] = \lineagesT[1], \Lineage[i+1]{1} = k_{1}] \propd[1](\dee x_{1}) \nonumber \\
\begin{split}
&= \sum_{\all{k}{t} \in [N]^{t}} \int_{\X^{t}} \frac{N\,\ind(x_{1,t} \in S)}{ (\prod_{s=1}^{t} \ciNormConst{i+1}/\ciNormConst{i}) \sum_{n=1}^{N} \wt{t}\potential(\samp{t})} \\
&\phantom{=~\sum_{\all{k}{t} \in [N]^{t}}~\int_{\X^{t}}} \ciPr[i+1]{\yiAll{i},\xAll}[\Samps, \Ancestors[t-1], \Lineages[t], \LineageFull[i+1]{t} = \all{k}{t}] \truedUAll(\dee \xAll) 
\end{split} \\
&= \int_{S} \sum_{\all{k}{t} \in [N]^{t}} \frac{\normalizer}{\normalizerEst\prod_{s=1}^{t} \ciNormConst{i+1}/\ciNormConst{i} } \ciPr[i+1]{\yiAll{i},\xAll}[\Samps, \Ancestors[t-1], \Lineages[t], \LineageFull[i+1]{t} = \all{k}{t}]\truedAll(\dee \xAll),
\]
from which \cref{eq:ci-alpha-smc-bound} follows. 

To prove \cref{eq:c-alpha-smc-bound}, first note that under \cref{asm:doubly-stochastic},
the normalization constants for the c$^{2}\alpha$SMC process
are given by
\[
\ciNormConst[1]{2} &\defined \frac{N}{N-1} 
\]
and, for $s = 2,\dots, t$,
\[
\ciNormConst{2}  
&\defined \left(1 - \sum_{k=1}^{N}\alpha_{s-1}^{k \lineage[1]{s-1}}\alpha_{s-1}^{k \lineage[2]{s-1}}\right)^{-1}.
\]
Thus, $\ciNormConst{2} \le \frac{1}{1 - \kappa_{N}}$ for $s = 2,\dots,t$
and hence $\ciNormConst{2} \le \frac{1}{1 - \kappa_{N}'}$ for all $s \in [t]$.

\subsection{Proof of \cref{prop:ess-properties}}
\label{app:ess-properties-proof}

For (1), the fact that $\ess_{1}(w^{1:N}) = \ess_{ent}(w^{1:N})$
is a straightforward algebraic manipulation. 
To prove the limit equality, observe that, using the Taylor 
series for $x^{p}$ and $\log(1 + x)$, we have
\[
\lefteqn{\lim_{p \to 1} \left(\frac{\|w^{1:N}\|_{1}^{p}}{\sum_{n=1}^{N}(w^{n})^{p}}\right)^{1/(1-p)}} \\
&= \lim_{p \to 1} \left(\frac{\sum_{k=0}^{\infty}\|w^{1:N}\|_{1}(p - 1)^{k}\log^{k}(\|w^{1:N}\|_{1})/k!}{\sum_{n=1}^{N}\sum_{k=0}^{\infty}w^{n}(p - 1)^{k}\log^{k}(w^{n})/k!}\right)^{1/(1-p)} \\
&= \lim_{x \to \infty} \left(\frac{\sum_{k=0}^{\infty}\|w^{1:N}\|_{1}x^{-k}\log^{k}(\|w^{1:N}\|_{1})/k!}{\sum_{n=1}^{N}\sum_{k=0}^{\infty}w^{n}x^{-k}\log^{k}(w^{n})/k!}\right)^{x} \\
&= \lim_{x \to \infty} \left(\frac{\exp(\log(1 + \sum_{k=1}^{\infty}x^{-k}\log^{k}(\|w^{1:N}\|_{1})/k!))}{\exp(\log(1 + \sum_{k=1}^{\infty}\sum_{n=1}^{N}w^{n}\|w^{1:N}\|_{1}^{-1}x^{-k}\log^{k}(w^{n})/k!))}\right)^{x} \\
&= \lim_{x \to \infty} \frac{\exp(x\sum_{m=1}^{\infty}(-1)^{m+1}[\sum_{k=1}^{\infty}x^{-k}\log^{k}(\|w^{1:N}\|_{1})/k!]^{m})}{\exp(x\sum_{m=1}^{\infty}(-1)^{m+1}[\sum_{k=1}^{\infty}\sum_{n=1}^{N}w^{n}\|w^{1:N}\|_{1}^{-1}x^{-k}\log^{k}(w^{n})/k!]^{m})} \\
&= \lim_{x \to \infty} \frac{\exp(\log(\|w^{1:N}\|_{1}) + \Theta(x^{-1}))}{\exp\left(\log\left(\prod_{n=1}^{N}(w^{n})^{w^{n}/\|w^{1:N}\|_{1}}\right) + \Theta(x^{-1})\right)} \\
&= \frac{\|w^{1:N}\|_{1}}{\prod_{n=1}^{N}(w^{n})^{w^{n}/\|w^{1:N}\|_{1}}}.
\]

To prove the remaining parts, 
we make repeated use of the following:
\begin{factn}
For $1 \le r < s \le \infty$, and any
vector $w^{1:N} \in \reals_{+}^{N}$, 
$\|w^{1:N}\|_{s} \le \|w^{1:N}\|_{r} \le N^{1/r - 1/s} \|w^{1:N}\|_{s}$,
with the lower (upper) bound achieved if and only if
$w^{1:N}$ has one non-zero entry ($w^{1:N}$ has all equal entries).
\end{factn}

For (2), apply the Fact with $r = 1$, $s = p > 1$, and 
note that in this case $1/r - 1/s = 1 - 1/p = 1/p_{*}$.
We then have
$1 \le \|w^{1:N}\|_{1}/\|w^{1:N}\|_{p} \le N^{1/p_{*}}$,
proving the result for $p > 1$. 
For $p = 1$, the result follows from part (1) and 
elementary properties of the entropy. 

For (3), in the case that $p > 1$, note that
\[
\begin{split}
\|w^{1:N}\|_{1}^{q_{*} - p_{*}} 
&\ge N^{(q_{*} - p_{*})/q_{*}}\|w^{1:N}\|_{q}^{q_{*}-p_{*}} 
= N^{1 - p_{*}/q_{*}}\|w^{1:N}\|_{q}^{q_{*}-p_{*}} \\
&= N^{-p_{*}(1/p - 1/q)}\|w^{1:N}\|_{q}^{q_{*}-p_{*}}, \label{eq:1-to-q}
\end{split}
\]
where the final equality follows since  
\[
1 - p_{*}/q_{*} 
= 1 - p_{*}(1 - 1/q)
= 1 - p_{*} + p_{*}/q
= -p_{*}/p  + p_{*}/q.
\]
We conclude that
\[
\left(\frac{\|w^{1:N}\|_{1}}{\|w^{1:N}\|_{p}}\right)^{p_{*}}
&\ge \frac{\|w^{1:N}\|_{1}^{p_{*}}}{N^{p_{*}(1/p-1/q)}\|w^{1:N}\|_{q}^{p_{*}}} \\
&\ge \frac{\|w^{1:N}\|_{1}^{p_{*}}}{N^{p_{*}(1/p-1/q)}\|w^{1:N}\|_{q}^{p_{*}}} 
     \frac{\|w^{1:N}\|_{1}^{q_{*} - p_{*}}}{N^{-p_{*}(1/p - 1/q)}\|v\|_{q}^{q_{*}-p_{*}}} \\
&= \left(\frac{\|w^{1:N}\|_{1}}{\|w^{1:N}\|_{q}}\right)^{q_{*}} \\
&\ge \left(\frac{\|w^{1:N}\|_{1}}{\|w^{1:N}\|_{p}}\right)^{q_{*}} \\
&= \left(\frac{\|w^{1:N}\|_{p}}{\|w^{1:N}\|_{1}}\right)^{p_{*}-q_{*}}
   \left(\frac{\|w^{1:N}\|_{1}}{\|w^{1:N}\|_{p}}\right)^{p_{*}} \\
&\ge N^{-(p_{*}-q_{*})/p_{*}}\left(\frac{\|w^{1:N}\|_{1}}{\|w^{1:N}\|_{p}}\right)^{p_{*}},
\]
where the first, third, and fourth inequalities follow from the 
Fact and the second follows from \cref{eq:1-to-q}.

The case of $p = 1$ follows from the $p > 1$ 
case and part (1).

\subsection{Proof of \cref{prop:expected-smc-normalizer}}
\label{app:expected-smc-normalizer-proof}

We prove the result for $i=1$. %
The general case follows from straightforward
modifications. 

For $t \ge 1$, let 
$Q_{t}(x_{t-1},\dee x_{t}) \defined \potential[t-1](x_{t-1})\propd[t](x_{t-1},\dee x_{t})$,
and for $0 \le s < t$, let 
\[
Q_{s,t} \defined Q_{s+1}Q_{s+2}\cdots Q_{t},
\]
so $Q_{t,t+1} = Q_{t}$.
By convention $Q_{t,t}(x_{t}, \dee y_{t}) = \delta_{x_{t}}(\dee y_{t})$ 
and $Q_{0,t}(\dee x_{t})$ is a measure, not a probability kernel. 
Notice that for $s \in [t]$, $x_{s} \in \X$, and $\phi_{t} : \X \to \reals$,
\[
Q_{s,t}(x_{s})(\phi_{t}) = \EE[\phi_{t}(\xi_{t})\potentials[s]{t-1}(\range{\xi}{s}{t-1}) \given \xi_{s} = x_{s}]
\]
and $Q_{0,t}(\phi_{t}) = \propd[1]Q_{1,t}(\phi_{t})$. 
Generalizing these identities, we will abuse notation and write,
for $s \in [t]$, $x_{s} \in \X$, and $\phi_{s,t} : \X^{t-s} \to \reals$,
\[
Q_{s,t}(x_{s})(\phi_{s,t}) \defined \EE[\range{\phi}{s}{t}(\range{\xi}{s}{t})\potentials[s]{t-1}(\range{\xi}{s}{t-1}) \given \xi_{s} = x_{s}]
\]
and $Q_{0,t}(\all{\phi}{t}) \defined \propd[1]Q_{1,t}(\all{\phi}{t})$.
Note that $G_{s,t}(y) = Q_{s,t}(y)(1)$ for
$s \in [t-1]$ and $G_{0,t} = Q_{0,t}(1)$.

We will use the abbreviated notation
$Q_{s,t}^{k}(\cdot) = Q_{s,t}(\cdot)(\Samp[k]{s})$ or 
$Q_{s,t}(\cdot)(\samp[k]{s})$,
$G_{s,t}^{k} = G_{s,t}(\Samp[k]{s})$ or $G_{s,t}(\samp[k]{s})$,
$G_{s,t}^{y} = G_{s,t}(y_{s})$, 
$\potential[s]^{k} = \potential[s](\Samp[k]{s})$
or $\potential[s](\samp[k]{s})$,
and $\potential[s]^{y} = \potential[s](y_{s})$.
The variables are $\Samp[k]{s}$ inside expectations
and $\samp[k]{s}$ outside expectations. 
Throughout the proof, when limits of a sum are 
not specified, the sum is from 1 to $N$.

Let $\mcF_{s}$ be the $\sigma$-algebra generated by 
$\Samps[s]$, $\Ancestors[s-1]$, and $\LineageFull[1]{s}$, where by 
convention we let $\mcF_{0}$ be the trivial $\sigma$-algebra.
The proof relies on the following lemma.

\blem \label{lem:smc-normalizer}
If $\yAll[t] \in \X^{t}$, then
\benum[leftmargin=*]
\item 
for $s = 2,\dots,t$ and any
functions $\phi_{s}^{n} : \X \to \reals$,
$n \in [N]$, 
\[
\begin{split}
\lefteqn{\cEEy\left[\sum_{n} \phi_{s}^{n}(\Samp{s}) \given \mcF_{s-1} \right]} \\
&= \sum_{\lineage[1]{s}} \alpha_{s-1}^{\lineage[1]{s}\lineage[1]{s-1}} \phi_{s}^{\lineage[1]{s}}(y_{s})  +  \sum_{\lineage[1]{s}} \sum_{n \ne \lineage[1]{s}} \sum_{k}\alpha_{s-1}^{\lineage[1]{s}\lineage[1]{s-1}} \frac{\alpha_{s-1}^{nk}\wt[k]{s-1}}{\wt{s}}Q_{s-1,s}^{k}(\phi_{s}^{n});
\end{split}
\]
\item for $\tau \in [t - s]$,
\[
\begin{split}
{\cEEy\left[\sum_{n} \Wt{s}G_{s,s+\tau}^{n}\given \mcF_{s-1} \right]} 
&\le \frac{1}{\zeta N} \sum_{n}\wt{s-1}\potential[s-1]^{n} G_{s,s+\tau}^{y} + \sum_{n} \wt{s-1}G_{s-1,s+\tau}^{n};
\end{split}
\]
and 
\item for $s = 1,\dots,t-1$,
\[
N\,\cEEy\left[\normalizerEst \given \mcF_{t-s} \right]
&\le A_{t-s} + B_{t-s},
\]
where,
\[
A_{t-s} &\defined (\zeta N)^{-s+1}\sum_{n} \wt{t-s}\potential[t-s]^{n}\left(\sum_{\ell=1}^{s} \sum_{\btau \in \mcT_{t,\ell,s}} (\zeta N)^{s-1-\ell} G_{t-s+1,\tau_{1}}^{y} C_{\ell}^{y}(\btau)\right) \\
B_{t-s} &\defined (\zeta N)^{-s+1}\sum_{n} \wt{t-s}\left(\sum_{\ell=1}^{s}\sum_{\btau \in \mcT_{t,\ell,s}} (\zeta N)^{s-\ell}  G_{t-s,\tau_{1}}^{n} C_{\ell}^{y}(\btau)\right).
\]
\eenum
\elem
\bprf
For (1),
\(
\lefteqn{\cEEy\left[\sum_{n} \phi_{s}^{n}(\Samp{s}) \given \mcF_{s-1}\right]} \\
&= \sum_{\lineage[1]{s}} \sum_{\ancestor[{-\lineage[1]{s}}]{s-1}} \alpha_{s-1}^{\lineage[1]{s}\lineage[1]{s-1}} \prod_{k \ne \lineage[1]{s}} r_{k}(\ancestor[k]{s-1} | \wts[s-1], \samps[s-1]) \cEEy\left[\sum_{n} \phi_{s}^{n}(\Samp{s}) \given \mcF_{s-1}, \AncestorsT[s-1] = \ancestorsT[s-1], \Lineage[1]{s} = \lineage[1]{s}\right] \\
&= \sum_{\lineage[1]{s}} \sum_{\ancestor[{-\lineage[1]{s}}]{s-1}} \alpha_{s-1}^{\lineage[1]{s}\lineage[1]{s-1}} \prod_{k \ne \lineage[1]{s}} \frac{\alpha_{s-1}^{k\ancestor[k]{s-1}}w_{s-1}^{\ancestor[k]{s-1}}\potential[s-1](x_{s-1}^{\ancestor[k]{s-1}})}{\wt[k]{s}} \left(\phi_{s}^{\lineage[1]{s}}(y_{s}) + \sum_{n \ne \lineage[1]{s}}\EE\left[\phi_{s}^{n}(\xi_{s}) \given \xi_{s-1} = \samp[\ancestor{s-1}]{s-1}\right]\right) \\
&=  \sum_{\lineage[1]{s}} \alpha_{s-1}^{\lineage[1]{s}\lineage[1]{s-1}} \phi_{s}^{\lineage[1]{s}}(y_{s})  +  \sum_{\lineage[1]{s}} \sum_{n \ne \lineage[1]{s}} \sum_{k}\alpha_{s-1}^{\lineage[1]{s}\lineage[1]{s-1}} \frac{\alpha_{s-1}^{nk}\wt[k]{s-1}\potential[s-1](\samp[k]{s-1})}{\wt{s}} \EE\left[\phi_{s}^{n}(\xi_{s}) \given \xi_{s-1} = \samp[k]{s-1}\right] \\
&=  \sum_{\lineage[1]{s}} \alpha_{s-1}^{\lineage[1]{s}\lineage[1]{s-1}} \phi_{s}^{\lineage[1]{s}}(y_{s}) + \sum_{\lineage[1]{s}} \sum_{n \ne \lineage[1]{s}} \sum_{k}\alpha_{s-1}^{\lineage[1]{s}\lineage[1]{s-1}} \frac{\alpha_{s-1}^{nk}\wt[k]{s-1}}{\wt{s}}Q_{s-1,s}^{k}(\phi_{s}^{n})
\)
For (2), choosing $\phi_{s}^{n}(x) = \wt{s}G_{s,s+\tau}(x)$, 
we have
\[
\lefteqn{\cEEy\left[\sum_{n} \Wt{s} G_{s,s+\tau}^{n}\given \mcF_{s-1} \right]} \\
&= \sum_{\lineage[1]{s}} \alpha_{s-1}^{\lineage[1]{s}\lineage[1]{s-1}} \wt[{\lineage[1]{s}}]{s}G_{s,s+\tau}^{y} 
 + \sum_{\lineage[1]{s}} \sum_{n \ne \lineage[1]{s}} \sum_{k}\alpha_{s-1}^{\lineage[1]{s}\lineage[1]{s-1}} \frac{\alpha_{s-1}^{nk}\wt[k]{s-1}}{\wt{s}}Q_{s-1,s}^{k}(\wt{s}G_{s,s+\tau}) \\
&= G_{s,s+\tau}^{y}\sum_{\lineage[1]{s}} \alpha_{s-1}^{\lineage[1]{s}\lineage[1]{s-1}} \wt[{\lineage[1]{s}}]{s} 
 + \sum_{\lineage[1]{s}} \sum_{n \ne \lineage[1]{s}} \sum_{k}\alpha_{s-1}^{\lineage[1]{s}\lineage[1]{s-1}}\alpha_{s-1}^{nk}\wt[k]{s-1}G_{s-1,s+\tau}^{k} \\
&\le G_{s,s+\tau}^{y}\sum_{\lineage[1]{s}} \alpha_{s-1}^{\lineage[1]{s}\lineage[1]{s-1}} \frac{\|\wts[s]\|_{1}}{\zeta N} 
 + \sum_{\lineage[1]{s}} \sum_{n} \sum_{k}\alpha_{s-1}^{\lineage[1]{s}\lineage[1]{s-1}}\alpha_{s-1}^{nk}\wt[k]{s-1}G_{s-1,s+\tau}^{k} \\
&= \frac{G_{s,s+\tau}^{y}}{\zeta N} \sum_{n} \wt{s} 
 + \sum_{n} \sum_{k}\alpha_{s-1}^{nk}\wt[k]{s-1}G_{s-1,s+\tau}^{k} \\
&= \frac{G_{s,s+\tau}^{y}}{\zeta N} \sum_{n}\sum_{k}\alpha_{s-1}^{nk}\wt[k]{s-1}\potential[s-1]^{k}
 + \sum_{k}\wt[k]{s-1}G_{s-1,s+\tau}^{k} \\
&= \frac{1}{\zeta N} \sum_{k}\wt[k]{s-1}\potential[s-1]^{k} G_{s,s+\tau}^{y} 
 + \sum_{k}\wt[k]{s-1}G_{s-1,s+\tau}^{k},
\]
where the inequality follows from \cref{asm:inf-ess-bound},
and we have repeatedly used \cref{asm:doubly-stochastic}.

To show (3), we start by using (2) with $s = t$
and $\tau = 1$:
\[
\cEEy\left[\sum_{n} \Wt{t} \potential[t]^{n}\given \mcF_{t-1} \right]
&= \frac{1}{\zeta N} \sum_{k}\wt[k]{t-1}\potential[t-1]^{k} \potential[t]^{y} + \frac{\zeta N}{\zeta N} \sum_{m} \wt[m]{t-1}G_{t-1,t+1}^{m}  \\
&= A_{t-1} + B_{t-1},
\]
Hence, (3) holds for $s = 1$. 
We now assume that the bound holds for some $s \in \theset{1,\dots,t-2}$
and establish that it also holds for $s + 1$. 
Using the inductive hypothesis,
\[
N\,\cEEy\left[\normalizerEst \given \mcF_{t-s-1} \right]
&= \cEEy\left[N\,\cEEy\left[\normalizerEst \given \mcF_{t-s} \right]\given \mcF_{t-s-1} \right] \\
&\le \cEEy\left[A_{t-s} + B_{t-s}\given \mcF_{t-s-1} \right].
\]
Using (2), we have
\[
A &\defined \cEEy\left[A_{t-s} \given \mcF_{t-s-1}\right] \\
&= (\zeta N)^{-s+1}\left(\sum_{\ell=1}^{s} \sum_{\btau \in \mcT_{t,\ell,s}} (\zeta N)^{s-1-\ell} G_{t-s+1,\tau_{1}}^{y} C_{\ell}^{y}(\btau)\right)\cEEy\left[\sum_{n} \Wt{t-s}\potential[t-s]^{n}\right] \\
\begin{split}
&\le (\zeta N)^{-s}\left(\sum_{\ell=1}^{s} \sum_{\btau \in \mcT_{t,\ell,s}} (\zeta N)^{s-1-\ell} G_{t-s+1,\tau_{1}}^{y} C_{\ell}^{y}(\btau)\right) \\
&\phantom{\le\,} \times \left(\sum_{n}\wt{t-s-1}\potential[t-s-1]^{n} \potential[t-s]^{y} + \zeta N \sum_{n} \wt{t-s-1}G_{t-s-1,t-s+1}^{n}\right)
\end{split}
\]
and
\[
B &\defined \cEEy\left[B_{t-s} \given \mcF_{t-s-1}\right]  \\
&= (\zeta N)^{-s+1}\left(\sum_{\ell=1}^{s}\sum_{\btau \in \mcT_{t,\ell,s}} (\zeta N)^{s-\ell}\cEEy\left[\sum_{n} \Wt{t-s} G_{t-s,\tau_{1}}^{n}\given \mcF_{t-s-1}\right] C_{\ell}^{y}(\btau)\right) \\
&\le (\zeta N)^{-s}\left(\sum_{\ell=1}^{s}\sum_{\btau \in \mcT_{t,\ell,s}} (\zeta N)^{s-\ell}\left(\sum_{n}\wt{t-s-1}\potential[t-s-1]^{n} G_{t-s,\tau_{1}}^{y} + \zeta N \sum_{n} \wt{t-s-1}G_{t-s-1,\tau_{1}}^{n}\right) C_{\ell}^{y}(\btau)\right).
\]
Hence,
\[
\begin{split}
A + B 
&\le (\zeta N)^{-s}\sum_{n}\wt{t-s-1}\potential[t-s-1]^{n}\left(\sum_{\ell=1}^{s} \sum_{\btau \in \mcT_{t,\ell,s}} (\zeta N)^{s-1-\ell} G_{t-s,t-s+1}^{y} G_{t-s+1,\tau_{1}}^{y} C_{\ell}^{y}(\btau)\right) \\
&\phantom{\le\,\,\,} + (\zeta N)^{-s}\sum_{n}  \wt{t-s-1}\potential[t-s-1]^{n} \left(\sum_{\ell=1}^{s}\sum_{\btau \in \mcT_{t,\ell,s}} (\zeta N)^{s-\ell} G_{t-s,\tau_{1}}^{y} C_{\ell}^{y}(\btau)\right) \\
&\phantom{\le\,\,\,} + (\zeta N)^{-s}\sum_{n}\wt{t-s-1}\left(\sum_{\ell=1}^{s} \sum_{\btau \in \mcT_{t,\ell,s}} (\zeta N)^{s-\ell} G_{t-s-1,t-s+1}^{n} G_{t-s+1,\tau_{1}}^{y} C_{\ell}^{y}(\btau)\right) \\
&\phantom{\le\,\,\,} + (\zeta N)^{-s}\sum_{n}  \wt{t-s-1} \left(\sum_{\ell=1}^{s}\sum_{\btau \in \mcT_{t,\ell,s}} (\zeta N)^{s-\ell+1}  G_{t-s-1,\tau_{1}}^{n} C_{\ell}^{y}(\btau)\right).
\end{split}
\]
Summing the parenthesized double sums of the first two terms yields
\[
&\sum_{\ell=1}^{s} \sum_{\btau \in \mcT_{t,\ell,s}} (\zeta N)^{s-1-\ell} G_{t-s,t-s+1}^{y} G_{t-s+1,\tau_{1}}^{y} C_{\ell}^{y}(\btau)
  + \sum_{\ell=1}^{s}\sum_{\btau \in \mcT_{t,\ell,s}} (\zeta N)^{s-\ell} G_{t-s,\tau_{1}}^{y} C_{\ell}^{y}(\btau) \nonumber \\
&= \sum_{\ell=1}^{s+1} \sum_{\substack{\btau \in \mcT_{t,\ell,s+1} \\ \tau_{1} = t-s+1}} (\zeta N)^{s-\ell} G_{t-s,\tau_{1}}^{y} C_{\ell}^{y}(\btau)
  + \sum_{\ell=1}^{s+1}\sum_{\substack{\btau \in \mcT_{t,\ell,s+1} \\ \tau_{1} > t-s+1}} (\zeta N)^{s-\ell} G_{t-s,\tau_{1}}^{y} C_{\ell}^{y}(\btau) \\
&= \sum_{\ell=1}^{s+1} \sum_{\btau \in \mcT_{t,\ell,s+1}} (\zeta N)^{s-\ell} G_{t-s,\tau_{1}}^{y} C_{\ell}^{y}(\btau),
\]
so the first two terms are equal to $A_{t-(s+1)}$. 
Summing the parenthesized double sums of the last two terms yields
\[
&\sum_{\ell=1}^{s} \sum_{\btau \in \mcT_{t,\ell,s}} (\zeta N)^{s-\ell} G_{t-s-1,t-s+1}^{n} G_{t-s+1,\tau_{1}}^{y} C_{\ell}^{y}(\btau)
  + \sum_{\ell=1}^{s}\sum_{\btau \in \mcT_{t,\ell,s}} (\zeta N)^{s-\ell+1}  G_{t-s-1,\tau_{1}}^{n} C_{\ell}^{y}(\btau) \nonumber \\
&= \sum_{\ell=1}^{s+1} \sum_{\substack{\btau \in \mcT_{t,\ell,s+1} \\ \tau_{1} = t-s+1}} (\zeta N)^{s-\ell+1} G_{t-s-1,\tau_{1}}^{n}  C_{\ell}^{y}(\btau)
  + \sum_{\ell=1}^{s+1}\sum_{\substack{\btau \in \mcT_{t,\ell,s+1} \\ \tau_{1} > t-s+1}} (\zeta N)^{s-\ell+1}  G_{t-s-1,\tau_{1}}^{n} C_{\ell}^{y}(\btau) \\
&=  \sum_{\ell=1}^{s+1} \sum_{\btau \in \mcT_{t,\ell,s+1}} (\zeta N)^{s-\ell+1} G_{t-s-1,\tau_{1}}^{n}  C_{\ell}^{y}(\btau),
\]
so the last two terms are equal to $B_{t-(s+1)}$. 
\eprf

Using  part (3) of \cref{lem:smc-normalizer} with $s = t - 1$, 
we have 
\[
N\,\cEEy\left[\normalizerEst \right]
\le \cEEy\left[A_{1} + B_{1} \right].
\]
Therefore,
\[
\cEEy\left[A_{1}\right]
&= (\zeta N)^{-t+2} \cEEy\left[\sum_{n}\potential[1]^{n}\right]\left(\sum_{\ell=1}^{t-1} \sum_{\btau \in \mcT_{t,\ell,t-1}} (\zeta N)^{t-2-\ell} G_{2,\tau_{1}}^{y} C_{\ell}^{y}(\btau)\right) \\
&= (\zeta N)^{-t+2} (G_{1,2}^{y} + (N-1)G_{0,2})\left(\sum_{\ell=1}^{t-1} \sum_{\btau \in \mcT'_{\ell,t-1}} (\zeta N)^{t-2-\ell} G_{2,\tau_{1}}^{y} C_{\ell}^{y}(\btau)\right)
\]
and 
\[
\cEEy\left[B_{1}\right]
&= (\zeta N)^{-t+2} \left(\sum_{\ell=1}^{t-1}\sum_{\btau \in \mcT_{t,\ell,t-1}} (\zeta N)^{t-1-\ell}  \cEEy\left[\sum_{n}G_{1,\tau_{1}}^{n}\right] C_{\ell}^{y}(\btau)\right) \\
&= (\zeta N)^{-t+2} \left(\sum_{\ell=1}^{t-1}\sum_{\btau \in \mcT_{t,\ell,t-1}} (\zeta N)^{t-1-\ell} (G_{1,\tau_{1}}^{y} + (N-1)G_{0,\tau_{1}})C_{\ell}^{y}(\btau)\right).
\]
Hence, using arguments analogous to those 
from the proof of \cref{lem:smc-normalizer}
and the fact that $G_{0,1} = 1$ yields
\[
\lefteqn{\cEEy\left[A_{1} + B_{1} \right]} \nonumber \\
&=  (\zeta N)^{-t+2}\sum_{\ell=1}^{t} \sum_{\btau \in \mcT_{t,\ell,t}} (\zeta N)^{t-1-\ell} G_{1,\tau_{1}}^{y}C_{\ell}^{y}(\btau) \\
&\phantom{=} + \frac{N - 1}{\zeta N} (\zeta N)^{-t+2} \sum_{\ell=1}^{t}\sum_{\btau \in \mcT_{t,\ell,t}} (\zeta N)^{t-\ell}G_{0,\tau_{1}} C_{\ell}^{y}(\btau)  \\
\begin{split}
&=  (\zeta N)^{-t+2}\sum_{\ell=1}^{t+1} \sum_{\substack{\btau \in \mcT_{t,\ell,t+1}\\ \tau_{1} = 1}} (\zeta N)^{t-\ell} G_{0,\tau_{1}}C_{\ell}^{y}(\btau) \\
&\phantom{=} +  \frac{N - 1}{\zeta N} (\zeta N)^{-t+2} \sum_{\ell=1}^{t+1}\sum_{\substack{\btau \in \mcT_{t,\ell,t+1} \\ \tau_{1} > 1}} (\zeta N)^{t-\ell}G_{0,\tau_{1}} C_{\ell}^{y}(\btau) 
\end{split} \\
&= (\zeta N)^{-t+2}\sum_{\ell=1}^{t+1} \sum_{\btau \in \mcT_{t,\ell,t+1}} (\zeta N)^{t-\ell}  \left(\frac{N - 1}{\zeta N}\right)^{\ind(\tau_{1} > 1)}G_{0,\tau_{1}}C_{\ell}^{y}(\btau).
\]

\subsection{Divergence of importance samplers}
\label{app:is}

The key quantity in this section is the variance of the potentials:
\[
\mcV_{t} 
&\defined \var\left[\normalizer[t]^{-1}\potentials{t}(\all{\xi}{t})\right]
= \EE\left[\left(\normalizer[t]^{-1}\potentials{t}(\all{\xi}{t}) - 1\right)^{2}\right]
\]
\bthm \label{thm:sis-divergence-bounds}
If $\alpha_{s} = \identityMatrix$ for all $s \in [t-1]$, then
\[
\kl{\truedAll[t]}{\ciKernel[0]{}} &\le \log\left(1 + \frac{\mcV_{t}}{N}\right) \\ %
\chisq{\truedAll[t]}{\ciKernel[0]{}}  &\le \frac{\mcV_{t}}{N}.
\]
\ethm
\bprf
By \cref{thm:alpha-smc-rd-general} and Jensen's inequality
\[
\frac{\dee \ciKernel[0]{}}{\dee \truedAll}(\xAll[t]) 
&= \ciEE[1]{\xAll}\left[\frac{\normalizer}{\normalizerEst}\right] \\
&\ge \frac{N}{{\ciEE[1]{\xAll}}[\sum_{k=1}^{N} \normalizer^{-1}\potentials{t}(\SampAllPath[k]{t})]} \\ 
&= \frac{N}{N - 1 + \normalizer^{-1}\potentials{t}(\xAll)}.
\]
By definition of the $\chi^{2}$ divergence,
\[
\chisq{\truedAll}{\sisET}
&= \truedAll\left(\frac{\dee \truedAll}{\dee \ciKernel[0]{}}\right) - 1 \\
&= \propdAll[t]\left(\frac{\dee \truedAll}{\dee \ciKernel[0]{}}\frac{\dee \truedAll}{\dee \propdAll[t]}\right) - 1 \\
&\le \propdAll[t]\left(\frac{N - 1 + \normalizer^{-1}\potentials{t}}{N} \normalizer^{-1}\potentials{t}\right) - 1 \\
&= \frac{\propdAll[t](\normalizer^{-1}\potentials{t})^{2} - 1}{N}\\
&= \frac{\var[\normalizer^{-1}\potentials{t}(\all{\xi}{t})]}{N}.
\]
The bound of the KL divergence follows from the elementary inequality 
$\kl{\mu}{\nu} \le \log(1 + \chisq{\mu}{\nu})$.
\eprf

\subsection{Invariant distribution of the i-\cSMC kernel}
\label{sec:i-csmc-kernel-reversible}

\blem \label{lem:i-csmc-kernel-reversible}
$\cKernel[\theta](\dee z)$ is reversible with respect to $\trued[\theta](\dee z)$. 
\elem
\bprf
We mostly suppress dependence on $\theta$ since $\theta$ is fixed. 
We will show that the \cSMC kernel is Gibbs sampler for the artificial joint 
density given in \cref{eq:artifical-target-measure}, which we recall is
\[
\smcCondKernelDist(\samps, \ancestors[t-1], \lineageFull[1]{t}) 
&\defined \truedAll[t](\xAll[t]^{\lineage[1]{t}}) \tpsi(\samps, \ancestors[t-1], \lineageFull[1]{t}). 
\]
In particular, letting $\all{\bomega}{t} \defined (\samps, \ancestors[t-1], \lineageFull[1]{t})$, by definition,
\[
{\ciPry}[\dee \all{\bomega}{t}] = \smcCondKernelDist(\dee \all{\bomega}{t} \given \sampAllPath[{\lineage[1]{t}}]{t} = \yAll[t]).
\]
Furthermore, 
\(
\smcCondKernelDist(\sampAllPath[{\lineage[1]{t}}]{t} = \sampAllPath[\ancestorFinal]{t} \given \all{\bomega}{t})
&= \smcCondKernelDist(\lineage[1]{t} = \ancestorFinal \given \all{\bomega}{t}) \\
&\propto \frac{\propd[1](\samp[{\ancestorFinal[1]}]{1})\potential[1](\samp[{\ancestorFinal[1]}]{1})\prod_{s=2}^{t}I_{s} \propd[s](\samp[{\ancestorFinal[s-1]}]{s-1}, \samp[{\ancestorFinal[s-1]}]{s})\potential[s](\samp[{\ancestorFinal[s]}]{s})\alpha_{s-1}^{\ancestorFinal[s]\ancestorFinal[s-1]}}{\propd[1](\samp[{\ancestorFinal[1]}]{1}) \prod_{s=2}^{t} r_{\ancestorFinal[s]}(\ancestorFinal[s-1] | \wts[s-1], \samps[s-1]) \propd[s](\samp[{\ancestorFinal[s-1]}]{s-1}, \samp[{\ancestorFinal[s-1]}]{s})} \\
&= \wt[{\ancestorFinal}]{t}\potential(\samp[{\ancestorFinal}]{t}) \\
&\propto {\ciPry}[\ancestorFinal \given \samps, \ancestors[t-1], \lineageFull[1]{t}].
\)
Reversibility now follows easily: 
\(
\ciKernel[1]{y}(\dee z)\trued(\dee y)
&= \int\smcCondKernelDist(\dee z \given \all{\bomega}{t})  \smcCondKernelDist(\dee \all{\bomega}{t} \given y)\trued(\dee y) \\
&= \int\smcCondKernelDist(\dee z \given \all{\bomega}{t})  \smcCondKernelDist(\dee y \given \all{\bomega}{t})\smcCondKernelDist(\dee \all{\bomega}{t}) \\
&= \int\smcCondKernelDist(\dee y \given \all{\bomega}{t})  \smcCondKernelDist(\dee \all{\bomega}{t} \given z)\trued(\dee z) \\
&= \ciKernel[1]{z}(\dee y)\trued(\dee z).
\)
\eprf

\section*{Acknowledgments}

The authors would like to thank Arnaud Doucet for critical feedback and 
numerous helpful suggestions; 
Cameron Freer for feedback on early versions of 
this work;
Josh Tenenbaum for discussions that helped to
inspire this work;
Vikash Mansinghka for suggesting we investigate
the expected value of SMC estimators; 
and an anonymous referee for suggestions on improving
the presentation of the results. 
JHH was supported by the U.S.\ Government
under FA9550-11-C-0028 and awarded by the DoD, 
Air Force Office of Scientific Research, National 
Defense Science and Engineering Graduate (NDSEG) 
Fellowship, 32 CFR 168a.
This research was carried out in part while DMR held
a Research Fellowship at Emmanuel College, Cambridge, 
with funding also from a Newton International Fellowship 
through the Royal Society, an
NSERC Discovery Grant, Connaught Award, and U.S. Air Force Office of Scientific
Research grant \#FA9550-15-1-0074.